\pdfoutput=1
\documentclass{article}

\usepackage{arxiv}
\usepackage{caption}
\usepackage{subcaption}
\usepackage{amsmath}
\usepackage{amssymb}
\usepackage{amsthm}
\usepackage{mathrsfs}
\usepackage{float}
\usepackage[export]{adjustbox}
\usepackage{hyperref}
\usepackage[numbers,sort&compress]{natbib}
\usepackage{xcolor}

\newcommand{\jump}[1]{[\![#1]\!]}

\newtheorem{proposition}{Proposition}

\title{Interface Energy and Phase Transformations: A Comparative Analysis of Cahn-Hilliard and CALPHAD-based Models in Ternary Substitutional Alloys}

\date{}
\author{Wolfgang Flachberger \\
	Department of Mathematics\\
	University of Innsbruck\\
	6020 Innsbruck, Austria \\
	\texttt{wolfgang.flachberger@uibk.ac.at} \\
    \And
	Thomas Antretter \\
	Institute of Mechanics\\
	Montanuniversit\"at Leoben\\
	8700 Leoben, Austria \\
	\texttt{thomas.antretter@unileoben.ac.at} \\
    \And
	Swaroop Gaddikere-Nagaraja \\
	Institute of Mechanics\\
	Montanuniversit\"at Leoben\\
	8700 Leoben, Austria \\
	\texttt{swaroop.gaddikere-nagaraja@unileoben.ac.at} \\
	\And
    Silvia Leitner \\
	Materials Center Leoben Forschung GmbH \\
	8700 Leoben, Austria \\
	\texttt{silvia.leitner@mcl.at}
    \And
    Manuel Petersmann \\
	Kompetenzzentrum Automobil- und \\
	Industrieelektronik GmbH \\
	Villach, Austria \\
	\texttt{manuel.petersmann@k-ai.at} \\
	\And
	Jiri Svoboda \\
	Institute of Physics of Materials \\
	Academy of Sciences of the Czech Republic \\
	616 00 Brno, Czechia \\
	\texttt{svobj@ipm.cz} \\
}

\renewcommand{\undertitle}{}

\begin{document}

\maketitle

\begin{abstract}
Diffusional phase transformations in alloys are modeled either by phase-field methods of Cahn--Hilliard type \cite{biner2017programming}, in which a double-well free energy and a gradient term impose an interface energy, or by reactive diffusion models \cite{svoboda2013new}, in which the driving forces follow from the CALPHAD methodology, the convex hull of the molar free energy, and the interface carries no energy. Their numerical schemes are tied to their free energies: phase-field schemes break down when the interface energy is set to zero, and reactive diffusion codes rely on finite differences. We derive, from the thermodynamic extremal principle, a single finite element scheme that covers both limits: the stationarity of the Rayleighian functional (free energy rate plus one half of the dissipation, with mass conservation imposed by multipliers) with respect to the rates, fluxes and multipliers at a frozen state. Fluxes and chemical affinities are primary unknowns, so no derivative of the affinity is required. The time-discrete scheme is the corresponding incremental minimization; it conserves mass exactly and satisfies a discrete energy--dissipation inequality for every convex molar free energy, including the convex hull without regularization. We make explicit that this last case is a degenerate parabolic problem of Stefan type, whose degeneracy can be removed, without altering the free energy, only by a moving mesh; since our aim is the coupling to continuum mechanics and damage on a fixed mesh, we operate at the boundary of degeneracy and show that mass conservation, the energy balance and the second law remain intact there. Thermodynamic consistency is proved for multi-component systems, a scheme relating Onsager and diffusion coefficients is proposed, and the role of interface energy is studied in binary and ternary examples, the latter with vacancies as a non-conserved component.
\end{abstract}


\section{Introduction}

\subsection{Motivation}
One of the challenges in modern materials science is developing a mathematical description of the evolution of a material's microstructure and morphology during production, thermal treatment, or usage, which can profoundly impact its mechanical properties and failure behavior. A significant portion of these processes can be attributed to diffusional phase transformations. Although mathematical descriptions are available, modeling and simulation still pose significant challenges in contemporary materials science due to the complex interactions and effects observed in relation to diffusion.
Moreover, many highly specialized models are available, but they are limited by their narrow focus on a specific phenomenon, which can make selecting a single model a compromise. One such example is the Cahn-Hilliard model \cite{cahn1958free}, which excels in simulating spinodal decomposition, with existing finite element methods available \cite{biner2017programming}. The Cahn-Hilliard model's most prominent feature is the free energy expressed in terms of a double-well potential, which drives phase separation by imposing interface energy. However, if this interface energy term is set to zero, the available finite element methods break down. On the other hand, there are methodologies, where the interface energy is implicitly assumed to be zero, such as those developed by Svoboda et al.\ \cite{svoboda2013new,svoboda2017incorporation}, where the CALPHAD method is used to derive the driving forces for diffusion, and the equations of linear thermodynamics are employed \cite{onsager1931reciprocal}. We will refer to these models as reactive diffusion models. These approaches neglect interface energy and often rely on simple finite difference methods, which are restricted to computations on regular grids.

Part of the reason why the two communities have developed separately is a matter of language. The Cahn--Hilliard equation is usually introduced with the textbook double well $\phi^2(1-\phi)^2$, whose minima lie at $\phi = 0$ and $\phi = 1$, so that the system strives for a complete separation of the components. For this choice the notions of mole fraction, phase fraction and order parameter coincide, and the free energy of the mixed state relative to the pure components, the ``bump'' of the double well, is at the same time the interface energy. In an alloy the minima of the molar free energy lie inside the interval $(0,1)$, at the compositions of the coexisting solid solutions, and the three notions are distinct: the mole fraction is the conserved quantity, the phase fraction follows from it by the lever rule, and an order parameter, if one is introduced to distinguish the phases, is a further, non-conserved variable \cite{miehe2014formulation}. Because the textbook double well blurs this distinction, it has contributed to the impression that reactive diffusion and phase-field models are different theories. They are not: as shown in Section~\ref{subsec:affinity}, both follow from the same free energy functional and the same variational principle, and reactive diffusion is the Cahn--Hilliard model with the convex hull of the molar free energy and without the gradient term, case (ii) introduced below. The double well is thus a model free energy specific to spinodal decomposition and not part of the method; the modelling decision is the choice of the molar free energy, for which the CALPHAD methodology provides the data.

A novel simulation scheme has been developed in this publication, capable of treating systems both with and without interface energy.
This new methodology enables the use of driving forces derived not only from double-well type potentials but also from the CALPHAD methodology, thereby allowing full control and direct utilization of thermodynamic data. To the authors' knowledge, an approach this versatile is not available in the literature.
The scheme is derived from the thermodynamic extremal principle (TEP) \cite{onsager1931reciprocal, ziegler1963some, svoboda2005application, fischer2014thermodynamic}: the evolution equations, and with them the discrete variational forms, follow from the stationarity of a single functional, the Rayleighian, whose extremal properties have recently been examined in \cite{flachberger2026variational}. We adopt the formulation and the nomenclature of that work throughout.
Additionally, the methodology enables systematic investigations into the influence of interface energy on phase morphologies and is designed to be broadly applicable to a wide range of phase transformation problems, including multi-phase, multi-component systems, with other rate-dependent phenomena such as vacancy generation and annihilation.
To demonstrate this, we consider a vacancy diffusion model \cite{svoboda2006diffusion} that accounts for vacancy sources and sinks. The crucial role of vacancies in facilitating diffusion in crystals is well-established \cite{manning1971correlation}, enabling atomic transport and phenomena such as the Kirkendall effect \cite{svoboda2017incorporation} and diffusion-controlled fatigue and void formation \cite{nishimura2004molecular}, particularly at high temperatures. In substitutional alloys, vacancies can be treated as an additional component, making this an ideal test case to showcase the methodology's capability for ternary systems.

\subsection{The Cahn-Hilliard model}

In their seminal work, Cahn and Hilliard (1958) \cite{cahn1958free}
presented a novel approach to define the Helmholtz energy of binary
alloys. They postulated that the total Helmholtz energy $\mathcal{F}[\phi]$ of a volume $V$ with nonuniform composition can be represented as a functional of the mole fraction $\phi$, comprising the molar free energy $f_0(\phi)$ of a uniform composition system and a gradient energy term $\tfrac{\kappa}{2} |\nabla \phi|^2$. This yields the characteristic free energy functional:

\begin{equation}
    \label{eq:Helmholtz}
    \mathcal{F}[\phi] = \frac{1}{\Omega} \int_{V}\left( f_0(\phi) + \frac{\kappa}{2} | \nabla \phi |^2\right) \; dV
\end{equation}

\begin{figure}
\centering
\includegraphics[max size={0.7\linewidth}{0.7\paperheight}]{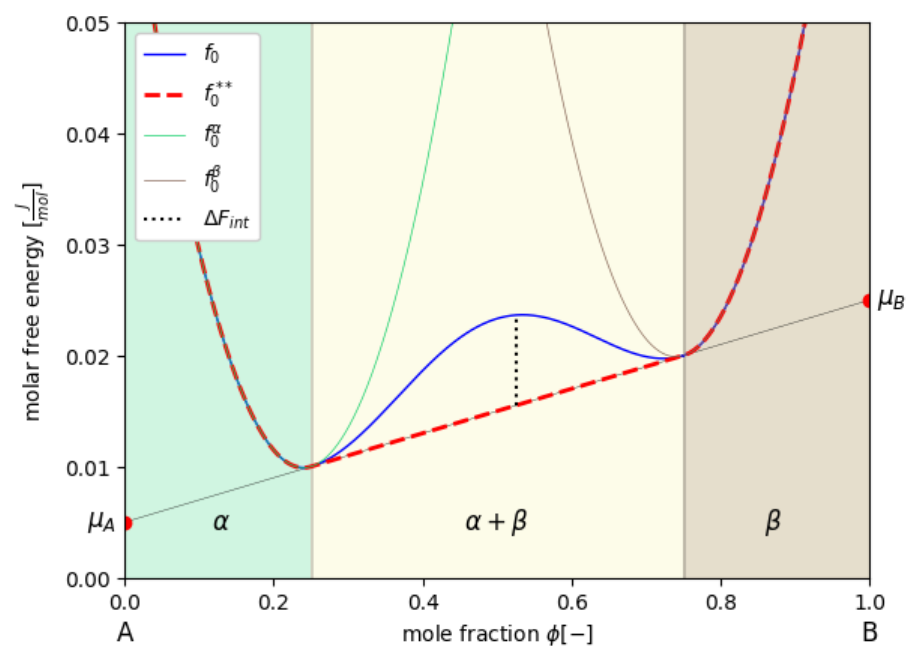}
\caption{ \(f_0\) represents a general double-well potential as
proposed by Cahn and Hilliard. The molar free energies of the individual
phases are given by \(f_{0}^{\alpha}\) and \(f_{0}^{\beta}\)
respectively. In an equilibrium system without interface energy the
molar free energy is given by the convex hull curve \(f_{0}^{**}\). Therefore, the maximal height difference $\Delta F_{int} $ between $f_0$ and $f_{0}^{**}$ represents the interface energy.}
\label{fig1}
\end{figure}

Here, \(\Omega\) denotes the molar volume of the alloy.
For a stress-free system with equal partial molar volumes of the components,
it can be considered constant.
For \(f_0(\phi)\), a double-well potential as given in Figure \ref{fig1} was proposed.
This also implies that the system accounts for interface energy, since the double-well suggests that
the mixed state (\(\phi=0.5\)) has a higher free energy than a linear
mixture of the individual phases in their equilibrium states
(\(\phi_{\alpha}=0.25\), \(\phi_{\beta}=0.75\)). If no interface energy is
intended, the double-well potential can be replaced by its convex hull
curve $f_{0}^{**}(\phi)$, as indicated by the dashed line in Figure \ref{fig1}.
This assumption is a widely adopted paradigm in physical chemistry and materials science, known as the common tangent construction or CALPHAD method \cite{lupis1983chemical}.
Using the convex
hull curve instead of the double-well potential introduces an intriguing
change in the behavior of a system's evolution, which will
no longer focus on minimizing the surface area
between the two phases. Instead, it will solely aim to minimize the
total free energy of the system through phase growth, subject to
boundary conditions and mass conservation.
We will refer to the replacement of $f_0(\phi)$ with $f_{0}^{**}(\phi)$ as case (i) in this publication.
Another noteworthy property of using the convex hull curve is that it can be employed without
regularization, i.e.\ $\kappa$ can be set to zero. The replacement of $f_0(\phi)$ with $f_{0}^{**}(\phi)$ together with $\kappa=0$ will be referred to as case (ii). Case (ii) results in a sharp transition between phases and has been extensively employed in various studies, including those by Svoboda et al.\ \cite{svoboda2006modeling}. The Cahn-Hilliard model, along with cases (i) and (ii), will be systematically compared in this work. The three cases are distinguishable solely by their respective molar free energy functions $f_m$, i.e., the integrand in equation \eqref{eq:Helmholtz}. Specifically, their molar free energies are defined as follows:

\begin{itemize}
    \item Cahn-Hilliard (interface energy + smooth interface)
        \begin{equation}
            \label{eq:case_ch}
            f_m(\phi) := f_0(\phi) + \tfrac{\kappa}{2} |\nabla \phi|^2
        \end{equation}
    \item Case (i): \emph{Reactive diffusion with regularization} (no interface energy + smooth interface)
        \begin{equation}
            \label{eq:case_1}
            f_m(\phi) := f_{0}^{**}(\phi) + \tfrac{\kappa}{2} |\nabla \phi|^2
        \end{equation}
    \item Case (ii): \emph{Reactive diffusion} (no interface energy + sharp interface)
        \begin{equation}
            \label{eq:case_2}
            f_m(\phi) := f_{0}^{**}(\phi)
        \end{equation}
\end{itemize}

Although Cahn and Hilliard referred to the term \(\tfrac{\kappa}{2} |\nabla \phi|^2\)
as ``gradient energy'', it has become common to call it a regularization
or homogenization term due to its smoothing effect on the interface.
Therefore, by neglecting regularization, the system tends towards discontinuous jumps in the mole fraction between the phases. As is made precise in Section~\ref{sec:degeneracy}, case (ii) is a degenerate parabolic problem: the diffusivity vanishes wherever the convex hull is linear, and the sharp interface is a free boundary. Schemes in which the mole fraction is the only unknown and the gradient of the affinity is required as a continuous field break down at this boundary, which is why numerical treatment by common diffusion schemes has up until now only been possible using the finite difference method (FD) \cite{svoboda2013new, svoboda2017incorporation} and a specific discontinuous Galerkin method which essentially mimics the previous FD schemes \cite{flachberger2024numerical}, both of which evaluate flux balances across cell faces instead. More elegant approaches to this problem are available \cite{dreyer2017sharp} but these methods are not within the scope of this work since we focus on methods compatible with the standard (continuous) finite element method (FEM), in order to be able to add mechanical coupling in future studies.
Interestingly, Cahn and Hilliard
themselves anticipated this case \cite{cahn1958free}:
\begin{quote}
\emph{We can see no reason why $\kappa$ should change sign or vanish in the unstable region, but we are unable to prove that it could not happen. If for some system there were a range over which $\kappa \leq 0$, then there would be a corresponding discontinuity in the interface profile and the treatment would have to be modified accordingly.}
\end{quote}
The ``unstable region'' is the spinodal region of the double well, in which a homogeneous composition is unstable and the two phases, and with them an interface, have to form. A vanishing $\kappa$ thus entails a discontinuous interface profile and calls for a modified treatment; Section~\ref{sec:degeneracy} makes precise what this modification amounts to for a convex molar free energy, and the scheme of Section~\ref{subsec:fem} is the modified treatment. For now, we
simply want to motivate the special case of \(\kappa = 0\) by noting
that it can be beneficial to have full control over the interface width
when performing multi-physics simulations, and that case (ii)
corresponds to employing the CALPHAD methodology. Moreover, it can be reasonable
to have a sharp interface when the length scale of the problem permits
it. Thermodynamics relies on the assumption that quantities like
pressure and temperature are assigned to a volume that is large enough
for the assumptions of statistical mechanics to hold. In the numerical
continuum treatment of non-equilibrium thermodynamics, this means that
each volume element in the domain must contain enough particles (atoms
or molecules) to ensure the validity of these assumptions, effectively
constraining the mesh size to be larger than a certain threshold.
However, when there is no interface energy and the length scale allows
it, it becomes legitimate to approximate the interface as sharp. Note
that this statement only holds for the continuum/volumetric treatment,
as there are also sharp interface models, such as those proposed by
\cite{dreyer2017sharp}, that operate with idealized surfaces representing
the interface. In these models, the interface is treated as a manifold,
allowing for the assignment of interface energy despite its
infinitesimally thin nature.

\subsection{Degeneracy of the sharp-interface limit and the scope of this work}
\label{sec:degeneracy}

We want to be explicit about the mathematical nature of case (ii), since it decides what can and what cannot be expected from any numerical method on a fixed mesh. As derived in Section~\ref{subsec:affinity}, the evolution of the mole fraction of a binary system is governed by $\frac{\partial \phi}{\partial t} = \nabla \cdot \big( M\, \nabla f_m'(\phi) \big)$ with a positive mobility $M$, so that the effective diffusivity is $M f_m''(\phi)$. For the convex hull $f_m = f_0^{**}$ this diffusivity vanishes identically on the two-phase interval $[\phi_\alpha, \phi_\beta]$, on which $f_0^{**}$ is linear: the equation is \textit{degenerate parabolic}. Its structure is that of the enthalpy formulation of the two-phase Stefan problem \cite{visintin1996models, elliott1987error, nochetto1988approximation}, $\frac{\partial \phi}{\partial t} = \Delta \beta(\phi)$ with a monotone $\beta = M f_0^{**\prime}$ that is constant on an interval: the mole fraction plays the part of the enthalpy, the affinity $-f_0^{**\prime}(\phi)$ that of the temperature (up to sign, an additive constant and a constant factor), and the width $\phi_\beta - \phi_\alpha$ of the two-phase interval that of the latent heat. Inside the two-phase region there is no driving force; the mole fraction there changes only through the divergence of the fluxes that the adjacent bulk phases impose at its boundary, and across the interface the mole fraction jumps from $\phi_\alpha$ to $\phi_\beta$ while the affinity is continuous. This degeneracy sits in the free energy and is distinct from the degenerate \textit{mobility} of the Cahn--Hilliard equation with $M = M(\phi)$ vanishing in the pure phases \cite{elliott1996cahn}, with which it should not be confused.

We are well aware of this degeneracy. It is not a defect of a particular discretization but a property of the continuum model. Regularization ($\kappa > 0$) removes it only by changing the free energy, at the price of an interface width that has no physical meaning for a convex free energy; without altering the free energy, it can be avoided only by abandoning the fixed-mesh, volumetric description, that is, by tracking the interface explicitly as a free boundary with a moving-mesh or front-tracking method \cite{baines2009moving, baines2011velocity}, in which the two-phase region has no volume and the mass balance across the interface is imposed as a jump condition. Such methods are not compatible with the scope of our work. Our aim is to couple diffusional phase transformations and damage to continuum mechanics, and a moving mesh would make the coupled problem too complicated. We therefore restrict ourselves to operating on the boundary of degeneracy: we accept the loss of regularity of the mole fraction at the phase boundary and ask only that the model conserves mass, closes the energy balance exactly, and enforces the second law of thermodynamics. Sections~\ref{sec:tep} and~\ref{subsec:fem} show that these three properties survive the degeneracy at the continuous and at the discrete level, for every convex molar free energy, the convex hull included: degeneracy at a phase boundary does not prevent us from enforcing the physics. This is sufficient for what we are doing. A better understanding of the interface itself would require atomistic models; we restrict ourselves to the thermodynamic, that is, the continuum description.

\section{Methodology}
\subsection{Affinity, chemical potential and thermodynamic consistency in the realm of a variational principle for multi-component systems}
\label{subsec:affinity}

Consider a thermodynamic system that is composed of \(n\)
distinct atomic components. The flux $\underline{j}_i$ and the chemical potential \(\mu_i\) of component $i$ are defined  for every $i\in\{1,\ldots,n\}$ as follows:

\begin{align}
    \label{eq:def_flux}
    \underline{j}_i &= - \sum_{k=1}^{n} \frac{L_{ik}}{T} \nabla \mu_k \\
    \label{eq:def_mu}
    \mu_{i} &= \Big( \frac{\partial F}{\partial m_i} \Big)_{T,V,m_{j \neq i}}
\end{align}

Here, \(m_i\) is the number of moles of component \(i\) in a reference
volume element and \(F\) denotes the total free
energy in the same volume element. \(L_{ik}\) and \(T\) represent the
Onsager coefficients and the temperature respectively. This is the usual
form in which the laws of diffusion are presented in the literature, including in \cite{kondepudi2014modern}.
Due to the second law of thermodynamics, the following condition must
hold at all times:

\begin{equation}
    \label{eq:def_xi}
    \xi := - \sum_{k=1}^{n} \underline{j}_k \cdot \nabla \big( \frac{\mu_k}{T} \big) \geq 0
\end{equation}

\(\xi\) is the rate of entropy production per unit volume. In a
non-equilibrium system it is greater than zero while it is equal to zero
in an equilibrium state. At constant temperature, $T\xi$ is the rate of energy dissipation per unit volume, and we denote its integral over the system by $\mathcal{P} = \int_V T \xi \, dV$, the \textit{dissipation}, which has units of power \cite{flachberger2026variational}. A system is said to be thermodynamically
consistent if the inequality \(\xi \geq 0\) is satisfied and if its
energy and mass are conserved. If a system is to be described by mole fractions
\(\phi_i\) instead of moles we can make use of the following relations:

\begin{align}
    \label{eq:total_number_moles}
    m &= \sum_{i=1}^{n} m_i \\
    \label{eq:mole_fraction}
    \phi_i &= \frac{m_i}{m}
\end{align}

For the mole fractions the following condition must hold:

\begin{equation}
    \label{eq:mole_frac_cond}
    \sum_{i=1}^{n} \phi_i = 1
\end{equation}

In this work, we focus on isothermal systems with constant molar volume, as in the Cahn-Hilliard model. We therefore impose:

\begin{equation}
    \Omega, T = \text{const.}
\end{equation}

Furthermore, we can make use of the following constraint:

\begin{equation}
\label{eq:flux_constraint}
\sum_{i=1}^{n}  \underline{j}_i = \underline{0}
\end{equation}

This equation is true in two cases: (I) For diffusion processes in
fluids and solids where all components have equal partial molar volumes
\cite{kondepudi2014modern}. The equation is then also known as ``no
volume flow'' constraint. (II) For diffusion of substitutional
components in crystalline solids, e.g.~substitutional alloys
\cite{manning1971correlation} (also for different partial molar volumes
of the components). The chemical
potential of a component in terms of the mole fractions is given by
\cite{lupis1983chemical}:

\begin{equation}
\label{eq:mu_i-x_i}
\mu_{i} =  f_m + \sum_{k=1}^{n-1} (\boldsymbol{\delta}_{ik} - \phi_k) \frac{\partial f_m}{\partial \phi_k},\quad  \forall i\in\{1,\ldots,n\}
\end{equation}

Here, \(f_m\) denotes a general molar free energy (convex or double-well or else) and $\boldsymbol{\delta}_{ik}$ is the Kronecker Delta. Note that equation
\eqref{eq:mu_i-x_i} also satisfies the Gibbs-Duhem equation \cite{svoboda2023note}. A proof can be found in
Appendix \ref{app:Gibbs-Duhem}. Mass conservation can be enforced through the continuity equation:

\begin{equation}
    \label{eq:continuity}
    \frac{\partial c_i}{\partial t} + \nabla \cdot \underline{j}_i = 0, \quad \forall i\in\{1,\ldots,n\}
\end{equation}

where the concentration $c_i$ of a component \(i\) is given by:

\begin{equation}
c_i = \frac{\phi_i}{\Omega}, \quad \forall i\in\{1,\ldots,n\}
\end{equation}

Since the system is defined in terms of mole fractions we can express
the \(n\)-th component in terms of the other \(n-1\) components by
rearranging equation \eqref{eq:mole_frac_cond}:

\begin{equation}
    \label{eq:x_n}
    \phi_n = 1 - \sum_{i=1}^{n-1} \phi_i
\end{equation}

Therefore, only \(n-1\) variables are needed to define the system
composition. For compactness, the collective mole fractions of the
different components, i.e.~the system composition, can also be expressed
as an array:

\begin{equation}
    \{\phi\} := \begin{bmatrix} \phi_{1} \\ . \\ . \\ \phi_{n-1} \end{bmatrix}
\end{equation}

The same applies to the flux of the \(n\)-th component if equation \eqref{eq:flux_constraint} is
modified:

\begin{equation}
    \label{eq:flux_constraint_2}
    \underline{j}_n = - \sum_{i=1}^{n-1}  \underline{j}_i
\end{equation}

By rearranging the definition of the volumetric rate of entropy
production \(\xi\), we can explicitly express the last element of the
sum over all \(n\) components

\begin{equation}
    \xi = - \frac{1}{T} \Big( \sum_{k=1}^{n-1} \underline{j}_k \cdot \nabla \mu_k + \underline{j}_n \cdot \nabla \mu_n \Big) \geq 0
\end{equation}

This gives rise to the opportunity of inserting the expression for the
flux of the \(n\)-th component:

\begin{equation}
    \xi = - \frac{1}{T} \Big( \sum_{k=1}^{n-1} \underline{j}_k \cdot \nabla \mu_k + \big( - \sum_{i=1}^{n-1} \underline{j}_i \big) \cdot \nabla \mu_n \Big) \geq 0
\end{equation}

Since the two sums are equal except for the indices of chemical
potentials the inequality can be simplified:

\begin{equation}
    \xi = - \frac{1}{T} \sum_{k=1}^{n-1} \underline{j}_k \cdot \big( \nabla \mu_k - \nabla \mu_n \big) \geq 0
\end{equation}

By reviewing equations \eqref{eq:def_flux} and \eqref{eq:def_mu} it becomes evident that the system is overdetermined due to the constraints \eqref{eq:mole_frac_cond} and \eqref{eq:flux_constraint}. The usual definition of fluxes according to equation \eqref{eq:def_flux} is therefore replaced:

\begin{equation}
\label{eq:new_ji}
    \underline{j}_i := - \sum_{k=1}^{n-1} \frac{\tilde{L}_{ik}}{T} \big( \nabla \mu_k - \nabla \mu_n \big), \quad \forall i\in\{1,\ldots,n-1\}
\end{equation}

This ensures the positivity of the rate of entropy production:

\begin{equation}
    \label{eq:xi_positive}
    \xi = \frac{1}{T^2} \Big( \sum_{i,k=1}^{n-1} \tilde{L}_{ik} \big( \nabla \mu_i - \nabla \mu_n \big) \cdot \big( \nabla \mu_k - \nabla \mu_n \big) \Big) \geq 0
\end{equation}

Note that the Onsager coefficient matrix was also changed from $L_{ik}$ of size $(n \times n)$ to $\tilde{L}_{ik}$ of size $(n-1 \times n-1)$. The definitions of flux and chemical potential in \eqref{eq:def_flux} and \eqref{eq:def_mu} were actually not altered, but rather adapted to the case of having to account for the constraints \eqref{eq:mole_frac_cond} and \eqref{eq:flux_constraint} which made it necessary to modify Onsager's coefficients. Both matrices, $L_{ik}$ and $\tilde{L}_{ik}$ however, must be positive definite and symmetric, by Onsager's reciprocal relations \cite{onsager1931reciprocal}. The result for the flux in equation \eqref{eq:new_ji} might appear arbitrary but things are clarified by investigating the term
\(\big( \nabla \mu_k - \nabla \mu_n \big)\). Since the chemical
potential of component \(i\) is defined in equation \eqref{eq:mu_i-x_i}, the chemical
potential of the last component \(n\) can be identified as:

\begin{equation}
    \label{eq:mu_n}
    \mu_{n} := f_m - \sum_{k=1}^{n-1} \phi_k \frac{\partial f_m}{\partial \phi_k}
\end{equation}

Combining equation \eqref{eq:mu_i-x_i} and \eqref{eq:mu_n} yields:

\begin{equation}
    \label{eq:aff}
    \mu_{i} - \mu_{n} = \frac{\partial f_m}{\partial \phi_i}, \quad \forall i\in\{1,\ldots, n-1\}
\end{equation}

Note that the negative of the partial derivative of a molar free energy quantity (Helmholtz or Gibbs) with respect to a non-dimensional variable which represents the extent of a reaction can be referred to as chemical affinity (as defined by the IUPAC, see also \cite{kondepudi2014modern}). The mole fraction, a dimensionless quantity, can be interpreted as the extent of a reaction, where values of zero or one correspond to full decomposition, representing the complete transformation from a mixed state to a decomposed state. Therefore, we denote the chemical affinity of the \(i\)-th component as \(\bar{\mu}_i\) (in accordance with equation \eqref{eq:aff}):

\begin{equation}
    \label{eq:Affinity}
    \bar{\mu}_i = - \frac{\partial f_m}{\partial \phi_i} = -(\mu_{i} - \mu_{n}) \quad \forall i\in\{1,\ldots,n-1\}
\end{equation}

The distinction between $\bar{\mu}_i$ and $\mu_i$ matters. The derivative of the free energy with respect to the mole fraction is routinely called ``chemical potential'' in the phase-field literature, but by \eqref{eq:aff} it is the difference of two chemical potentials: the driving force for the exchange of an atom of component $i$ against an atom of component $n$, known as the diffusion potential in the work of Larch\'e and Cahn \cite{larche1985interactions}, and an affinity in the sense of the IUPAC definition quoted above. Identifying $\bar{\mu}_i$ with $\mu_i$ makes the parameters of a phase-field model hard to relate to measured chemical potentials and to CALPHAD data, and it is one source of the recurring debate on whether such models are thermodynamically consistent or merely phenomenological. We therefore keep the notation of materials science, in which $\mu_i$ is the chemical potential of \eqref{eq:def_mu}, and reserve $\bar{\mu}_i$ for the affinity.

These derivations are also possible if a regularization term is included in $f_m$, which simply requires replacing all partial derivatives with variational derivatives, i.e., applying the Euler-Lagrange equations. Equation \eqref{eq:Affinity} can therefore be replaced by:

\begin{equation}
    \label{eq:Affinityvariational}
    \bar{\mu}_i = - \frac{\delta f_m}{\delta \phi_i} = -(\mu_{i} - \mu_{n}) \quad \forall i\in\{1,\ldots,n-1\}
\end{equation}

Combining equations \eqref{eq:Affinityvariational} and \eqref{eq:new_ji} further simplifies the definition of the fluxes:

\begin{equation}
    \label{eq:derived_j_}
    \underline{j}_k = - \sum_{i=1}^{n-1} \frac{\tilde{L}_{ik}}{T} \nabla \big(  \frac{\delta f_m}{\delta \phi_i} \big), \quad \forall k\in\{1,\ldots,n-1\}
\end{equation}

It is thus evident, as noted earlier, that reactive diffusion models and the Cahn-Hilliard model are fundamentally equivalent, just distinguishable by their definitions of the molar free energy
\(f_m\) (see equations \eqref{eq:case_ch} - \eqref{eq:case_2}). Considering a binary system and combining equations \eqref{eq:continuity} and \eqref{eq:derived_j_} yields:

\begin{equation}
    \label{eq:cahn-hilliard}
    \frac{\partial c}{\partial t} = \nabla \cdot \Big( \frac{\tilde{L}}{T} \nabla \Big( \frac{\delta f_m}{\delta \phi} \Big) \Big)
\end{equation}

Equation \eqref{eq:cahn-hilliard} closely resembles the usual way the Cahn-Hilliard model is presented
in the literature \cite{biner2017programming, nagaraja2023variational}. However, it should be noted that the parameters $\tilde{L}, T$ and $\Omega$ are often unified into a single mobility parameter to obtain a simple evolution equation for $\phi$. Unfortunately, this can lead to confusion regarding units, which is why we believe that the derivation and proof of thermodynamic consistency in this work were necessary. While proofs of thermodynamic consistency are also utilized in other works, such as \cite{emmerich2003diffuse}, their derivations are highly formalized, omitting units and restricted to ``extensive, locally conserved variables'' instead of a specific quantity, such as the mole fraction $\phi$. The derivations herein therefore close a gap and also shed light on the role of the quantity $\bar{\mu}$ as defined in equation \eqref{eq:Affinityvariational}, which is often mistaken for the chemical potential. If the mole fraction defines the system composition and constraints of type \eqref{eq:flux_constraint} are considered, the variation of the molar free energy \eqref{eq:Affinityvariational} should be referred to as affinity. However, specific problem formulations may exist where this quantity indeed represents a chemical potential.

Equation \eqref{eq:cahn-hilliard} also exhibits the degeneracy of Section~\ref{sec:degeneracy}. In its primal form the equation requires the gradient of the affinity, $\nabla (\delta f_m / \delta \phi) = f_m''(\phi) \nabla \phi$ for $\kappa = 0$, so that the differential operators demand $C^3$ continuity of $f_m$ in the classical sense and $C^2$ continuity in a variational formulation\footnote{This is also necessary since the convex hull curve $f_{0}^{**}$ is only $C^1$-continuous: its second derivative jumps at $\phi_\alpha$ and $\phi_\beta$ and vanishes in between.}; moreover, $f_m''$ vanishes on the two-phase interval of the convex hull, where the equation ceases to be parabolic and determines $\phi$ only through mass conservation. Across the interface, $\phi$ is discontinuous while $\delta f_m / \delta \phi$ is continuous. Any scheme that requires $\nabla (\delta f_m / \delta \phi)$ as a continuous field therefore fails for case (ii). Regularization ($\kappa > 0$, case (i)) restores the parabolic character and the regularity of $\phi$, at the price of an interface width that has no physical meaning for a convex free energy; the scheme of Section~\ref{subsec:fem} avoids the differentiation of the affinity altogether and treats both cases alike.

\subsection{The thermodynamic extremal principle and the Rayleighian}
\label{sec:tep}

All equations defining the diffusion problem, \eqref{eq:continuity}, \eqref{eq:Affinityvariational} and \eqref{eq:derived_j_}, follow from a single variational statement: the thermodynamic extremal principle (TEP). The idea goes back to Onsager's principle of least dissipation of energy \cite{onsager1931reciprocal} and Ziegler's principle of maximum rate of dissipation \cite{ziegler1963some}, was related to dissipation potentials by Hackl and Fischer \cite{hackl2008relation}, and has been systematized for materials science by Svoboda et al.\ \cite{svoboda2005application, fischer2014thermodynamic}; see also \cite{emmerich2003diffuse} for phase-field models and \cite{doi2011onsager, beris2024dissipation} for the same principle in soft matter and fluid mechanics, where it is known as the Onsager variational principle. We follow the formulation of \cite{flachberger2026variational}, where the Rayleighian is derived from the constrained stationarity of the dissipation and the extremal properties of the principle are examined, and use its nomenclature. That work is concerned with fluids and mechanics, where inertia enters, and leaves diffusion aside; the present work applies the same formulation to multi-component diffusion.

Consider a closed, isothermal system with the free energy $\mathcal{F} = \frac{1}{\Omega}\int_V f_m \, dV$ of equation \eqref{eq:Helmholtz}. The TEP requires the rate of the free energy. For a crystalline solid whose lattice is at rest the control volume is fixed in space, and the Reynolds transport theorem reduces to the chain rule:
\begin{equation}
    \label{eq:F_rate}
    \dot{\mathcal{F}} = \frac{1}{\Omega} \int_V \sum_{k=1}^{n-1} \frac{\delta f_m}{\delta \phi_k}\, \frac{\partial \phi_k}{\partial t} \; dV = -\frac{1}{\Omega} \int_V \sum_{k=1}^{n-1} \bar{\mu}_k\, \frac{\partial \phi_k}{\partial t} \; dV\, .
\end{equation}
The rate is linear in the rates $\partial \phi_k / \partial t$, and the coefficients, the driving forces $\delta f_m / \delta \phi_k = -\bar{\mu}_k$, are state functions: they are evaluated at the given state $\{\phi\}$ and are not varied. (For $\kappa > 0$ the variational derivative in \eqref{eq:F_rate} arises from the chain rule after an integration by parts with $\nabla \phi_k \cdot \underline{n} = 0$ on $\partial V$; Section~\ref{subsec:fem} uses the directional form, which needs no boundary condition.) The dissipation $\mathcal{P} = \int_V T\xi \, dV = T \int_V \sum_{i,k=1}^{n-1} \tilde{L}^{-1}_{ik}\, \underline{j}_i \cdot \underline{j}_k \; dV$ is a positive-definite quadratic form in the fluxes, with $\tilde{L}^{-1}_{ik}$ the kinetic resistance, and the energy balance of the closed system reads $\dot{\mathcal{F}} + \mathcal{P} = 0$; the second law enters through $\mathcal{P} \geq 0$. The TEP selects, among all fluxes compatible with this balance, those that render the dissipation stationary. As shown in \cite{flachberger2026variational} for generalized rates, of which the fluxes are an instance once the rates $\partial \phi_k / \partial t$ are expressed through them by mass conservation, and by Svoboda, Hackl and Fischer \cite{svoboda2023note} for diffusion, the balance defines an ellipsoid through the origin in the space of the fluxes, and the physical fluxes \textit{maximize} the dissipation on it (Ziegler's principle); equivalently, they are the unconstrained minimizer of the convex functional
\begin{equation}
    \label{eq:rayleighian_def}
    \mathcal{R} = \dot{\mathcal{F}} + \tfrac{1}{2} \mathcal{P}\, ,
\end{equation}
the \textit{Rayleighian functional} \cite{rayleigh1873some, doi2011onsager, flachberger2026variational}. Mass conservation \eqref{eq:continuity} is a constraint on the rates and fluxes, and it is added to the Rayleighian with the multipliers $\bar{\mu}_k$. The Rayleighian \textit{density} of the multi-component diffusion problem thus reads:

\begin{equation}
    \label{eq:rayleighian}
    \mathscr{R}(\{\frac{\partial \phi}{\partial t}\},\{\underline{j}\},\{\bar{\mu}\}) := \frac{1}{\Omega} \sum_{k=1}^{n-1} \frac{\delta f_m}{\delta \phi_k} \, \frac{\partial \phi_k}{\partial t} + \frac{T}{2} \sum_{i,k=1}^{n-1} \tilde{L}_{ik}^{-1} (\underline{j}_i \cdot \underline{j}_k) + \sum_{k=1}^{n-1} \bar{\mu}_k \Big(  \frac{1}{\Omega} \frac{\partial \phi_k}{\partial t} + \nabla \cdot \underline{j}_k \Big)
\end{equation}

with $\mathcal{R} = \int_V \mathscr{R} \, dV$. Following \cite{flachberger2026variational}, script letters denote densities and calligraphic letters the functionals of the whole system. The first term of \eqref{eq:rayleighian} is the rate of the free energy \eqref{eq:F_rate}, the second is one half of the dissipation, and the third is the constraint of mass conservation. The same density was called a Lagrangian in our earlier work \cite{flachberger2024numerical, nagaraja2023variational, flachberger5251955impact}; we adopt the name Rayleighian here because the TEP is fundamentally different from Hamilton's principle, despite the shared vocabulary of Euler--Lagrange equations: it varies the rates at a single instant, for a given state, and its stationarity conditions are algebraic or spatial relations between the rates, with no temporal differentiation \cite{flachberger2026variational}. Four rules make the construction unambiguous \cite{flachberger2026variational}: (a) the energy rate is fixed by the Reynolds transport theorem with the control volume appropriate to the physics, here the fixed volume of a lattice at rest; (b) only the rates $\partial \phi_k / \partial t$, the fluxes $\underline{j}_k$ and the multipliers $\bar{\mu}_k$ are varied, at the current instant, and no temporal Euler--Lagrange equations are applied; (c) the state $\{\phi\}$ is given, so that the driving forces and the Onsager coefficients are evaluated at it and are not varied; (d) driving forces that contain accelerations are held fixed during the variation. Rule (d) is not needed here: diffusion is inertialess, the driving forces are state functions, and nothing has to be frozen by hand, which is the situation in which the TEP is most elegant.

Stationarity of $\mathcal{R}$ with respect to $\partial \phi_i / \partial t$, $\underline{j}_i$ and $\bar{\mu}_i$ yields equations \eqref{eq:Affinityvariational}, \eqref{eq:derived_j_} and \eqref{eq:continuity}, respectively: the first identifies the multiplier with the affinity, the second is Onsager's law, in which the gradient of the multiplier appears after an integration by parts of the constraint term, and the third is the continuity equation. The stationary point is a minimum with respect to the fluxes (and, through the constraint, the rates) and a maximum with respect to the multipliers, that is, a saddle point of the augmented functional, as for every Lagrange multiplier. The Rayleighian \eqref{eq:rayleighian} is strictly convex in the fluxes and linear in the rates; the minimization of the Rayleighian at each instant is the time-continuous form of the minimizing-movement scheme by which gradient flows are constructed \cite{jordan1998variational, otto2001geometry, flachberger2026variational}, a fact that we exploit in Section~\ref{subsec:fem}.
A detailed derivation of the principle, as presented in equation \eqref{eq:rayleighian}, can be found in the works of Svoboda et al.\ \cite{svoboda2005application} or by Emmerich \cite{emmerich2003diffuse}.
The concept is also particularly useful for exploring coupled problems of reactive diffusion and continuum mechanics \cite{flachberger5251955impact}: the Rayleighian of a coupled problem is obtained by summing the energy rates, the dissipation potentials and the constraints of the individual phenomena \cite{flachberger2026variational}. In this publication, we demonstrate how the Rayleighian can be used to automatically derive variational formulations for novel finite element schemes in reactive diffusion and phase-field modeling.

\subsection{Correlation of diffusion coefficients and Onsager's coefficients}
A question that often arises and is of great interest when applying phase-field methods or reactive diffusion methods is to find the dependence of Onsager's coefficients $\tilde{L}_{ik}$ on experimentally measured diffusion coefficients $D_i$ for all $i,k\in\{1,\ldots,n\}$. This is, however, not a trivial task and can hardly be answered in the general case. Fortunately, the investigated cases and the introduced methodology offer the opportunity to formulate $f_m$ in a specific way that allows straightforward relations between $\tilde{L}_{ik}$ and $D_i$. This is because both a double well and its convex counterpart can be constructed from quadratic functions that represent the free energies of individual phases. For example, the molar free energy of phase $\alpha$ can be represented as (see also Figure \ref{fig1}):

\begin{equation}
    f_{0}^{\alpha}(\phi) := f_{\alpha} + \frac{k}{2} (\phi-\phi_{eq}^{\alpha})^2
\end{equation}

Here, $f_{\alpha}$ represents the molar free energy of phase $\alpha$ in equilibrium composition (that is, $\phi=\phi_{eq}^{\alpha}$) and $k$ controls the curvature of its convex energy. In general $f_{\alpha}$ is dependent on temperature and volume ($T,V$) and can be calculated from heat capacity and mechanical properties of a material, but will be considered constant in the present work. This simple quadratic dependency is necessary to give a straightforward dependency on the diffusion coefficients since measuring diffusion coefficients often utilizes Fick's law of diffusion \cite{kumar2011intrinsic}. Therefore, it is required that Fick's law and the thermodynamic definition of flux lead to the same result. For a binary system this means, with $\bar{\mu} = -\partial f_m / \partial \phi$ according to \eqref{eq:Affinity}:

\begin{equation}
    \underline{j} := \frac{\tilde{L}}{T} \nabla \bar{\mu} = - \frac{\tilde{L}}{T} \nabla \Big( \frac{\partial f_m}{\partial \phi} \Big) = - \tilde{D} \nabla c
\end{equation}

Here $\tilde{D}$ denotes the inter-diffusion coefficient of components 1 and 2. The relation further simplifies to:

\begin{equation}
    \frac{\tilde{L}}{T} \nabla \Big( \frac{\partial f_m}{\partial \phi} \Big) = \tilde{D} \nabla \Big( \frac{\phi}{\Omega} \Big)
\end{equation}

If $f_m$ is constructed from quadratic functions, the following simplification holds at least in the different bulk phases but not at interfaces:

\begin{equation}
    \frac{\tilde{L}}{T} \nabla \Big( k \cdot (\phi-\phi_{eq}) \Big) = \tilde{D} \nabla \Big( \frac{\phi}{\Omega} \Big)
\end{equation}

Since applying derivatives to constants yields zero, a straightforward relation can be found:

\begin{align}
    \label{eq:onsager-diff}
    \frac{\tilde{L}k}{T} \nabla \phi &= \frac{\tilde{D}}{\Omega} \nabla \phi \\
    \tilde{L} &= \frac{\tilde{D}T}{k\Omega}
\end{align}

Obviously, the relation can be criticized for relating the Onsager coefficient $\tilde{L}$ correctly only in the bulk phases but not at interfaces, where $f_m$ is a linear function for a convex potential, and not a quadratic function. However, during simulation, the domain is dominated by bulk phases, and the interfaces only amount to a small fraction, which underscores the potential of this simplification. This is especially true for case (ii), where the representation of the interface is very sharp. Moreover, interface regions in materials science are known for their increased diffusional activity, which could be argued to be captured by this methodology. Additionally, it is worth noting the influence of the parameter $k$. Its role in the construction of the molar free energy $f_m$ is clear, but since it appears in the Onsager coefficient, it effectively loses its impact in determining the overall kinetics of the model. It can therefore be viewed as a phenomenological parameter that should be chosen to obtain good convergence of the simulation but has no direct consequence for the system's evolution.

\subsection{Scaling the Rayleighian}
\label{sec:scaling}
In the context of PDEs and simulation, \emph{scaling} refers to the process of expressing a physical problem in a canonical form \cite{langtangen2016scaling}. This may be achieved by viewing the equations in a non-dimensional setting, changing the unit of field variables and/or formulating the problem with the least possible number of physical parameters. Since a variational description of the system was chosen, the change of dimension is applied at the same fundamental level. Consider the Rayleighian of a binary system:

\begin{equation}
    \mathscr{R}(\frac{\partial \phi}{\partial t}, \underline{j}, \bar{\mu}) := \frac{1}{\Omega}  \frac{\delta f_m}{\delta \phi} \, \frac{\partial \phi}{\partial t} + \frac{T}{2\tilde{L}} |\underline{j}|^2 + \bar{\mu} \Big(\frac{1}{\Omega} \frac{\partial \phi}{\partial t} + \nabla \cdot \underline{j}\Big)
\end{equation}

Utilizing the derived dependency of the Onsager coefficient from equation \eqref{eq:onsager-diff} yields the following expression:

\begin{equation}
    \mathscr{R}(\frac{\partial \phi}{\partial t}, \underline{j}, \bar{\mu}) := \frac{1}{\Omega}  \frac{\delta f_m}{\delta \phi} \, \frac{\partial \phi}{\partial t} + \frac{k\Omega}{2\tilde{D}} |\underline{j}|^2 + \bar{\mu} \Big( \frac{1}{\Omega} \frac{\partial \phi}{\partial t} + \nabla \cdot \underline{j}\Big)
\end{equation}

By reviewing the equations that are obtained by stationarity with respect to $\underline{j}$ and $\bar{\mu}$, it is evident that the Rayleighian holds potential for simplification:

\begin{align}
    \label{eq:first_order_system_1}
    \underline{j} &= \frac{\tilde{D}}{k \Omega} \nabla \bar{\mu} \\
    \label{eq:first_order_system_2}
    \frac{1}{\Omega}\frac{\partial \phi}{\partial t} &= - \nabla \cdot \underline{j}
\end{align}

If equation \eqref{eq:first_order_system_1} and \eqref{eq:first_order_system_2} are combined into a second order equation, it can be seen that if the molar volume $\Omega$ is constant, it can be dropped from the equations entirely:

\begin{equation}
    \label{eq:second_order_system}
    \frac{1}{\Omega}\frac{\partial \phi}{\partial t} = - \nabla \cdot \Big( \frac{\tilde{D}}{k \Omega} \nabla \bar{\mu} \Big) = \nabla \cdot \Big( \frac{\tilde{D}}{k \Omega} \nabla \frac{\delta f_m}{\delta \phi} \Big)
\end{equation}

The same applies for the variation of the Rayleighian with respect to $\partial \phi / \partial t$. To this end the scaled Rayleighian for a binary system is introduced:

\begin{equation}
    \label{eq:scaled}
    \mathscr{R}(\frac{\partial \phi}{\partial t}, \underline{j}, \bar{\mu}) := \frac{\delta f_m}{\delta \phi} \, \frac{\partial \phi}{\partial t} + \frac{k}{2\tilde{D}} |\underline{j}|^2 +  \bar{\mu} \Big( \frac{\partial \phi}{\partial t}  + \nabla \cdot \underline{j}\Big)
\end{equation}

Note that the unit of flux is changed in the scaled Rayleighian from $[\frac{mol}{m^2s}]$, which corresponds to the usual definition, to $[\frac{m}{s}]$. Therefore the vector field $\underline{j}$ could also be referred to as velocity or drift velocity of the atomic components in the scaled Rayleighian. This ensures that the mole fraction and the chemical affinity can remain in their respective standard unit. In the scaled form, the dissipation is $\mathcal{P} = \int_V \frac{k}{\tilde{D}} |\underline{j}|^2 \, dV$, and the second term of \eqref{eq:scaled} is one half of its density, as in \eqref{eq:rayleighian_def}.

\subsection{Time discretization and variational formulation}
\label{subsec:fem}

For simulation purposes the Rayleighian \eqref{eq:rayleighian}, or its scaled counterpart \eqref{eq:scaled}, is adapted to a discrete time setting. We follow the rules of Section~\ref{sec:tep} literally: the variation is carried out first, at the frozen state, and the time discretization is applied afterwards. This removes an ambiguity of the earlier version of this work, in which a total time difference of the free energy was varied with respect to the new state; the resulting discrete forms are the same (see below), but the derivation is now unambiguous. To illustrate the process in a straightforward manner, we consider a binary system, and we write $f \in \{f_0, f_0^{**}\}$ for the bulk part of the molar free energy, so that $f_m = f(\phi) + \tfrac{\kappa}{2}|\nabla \phi|^2$, and $\mathcal{F}[\phi] = \int_V f_m \, dV$ in the scaled units of \eqref{eq:scaled}. At a given state $\phi$ the scaled Rayleighian functional reads

\begin{equation}
    \label{eq:rayleighian_functional_binary}
    \mathcal{R}\Big[\frac{\partial \phi}{\partial t}, \underline{j}, \bar{\mu}\Big] = \int_V \Big( \frac{\partial f}{\partial \phi}\, \frac{\partial \phi}{\partial t} + \kappa \nabla \phi \cdot \nabla \frac{\partial \phi}{\partial t} + \frac{k}{2\tilde{D}} |\underline{j}|^2 +  \bar{\mu} \Big( \frac{\partial \phi}{\partial t}  + \nabla \cdot \underline{j}\Big) \Big) \; dV,
\end{equation}

where the first two terms are the rate of the free energy, $\int_V (\delta f_m / \delta \phi)\, (\partial \phi / \partial t)\, dV$, written as the directional derivative of $\mathcal{F}$ in the direction of the rate, so that no integration by parts is needed and the gradient term can be included in the Cahn--Hilliard case and in case (i). The state $\phi$ is given. The unknowns are the rate $\partial \phi / \partial t$, the flux $\underline{j}$ and the multiplier $\bar{\mu}$, and the variation is carried out with respect to these three fields only, at the current instant. Denoting the directions of the variation by $\hat{\phi}$ for the rate, $\underline{\hat{j}}$ for the flux and $\hat{\bar{\mu}}$ for the multiplier, the stationarity conditions read

\begin{align}
    \label{eq:stat_1}
    0 &= \int_{V} \Big( \frac{\partial f(\phi)}{\partial \phi} \hat{\phi} + \kappa \nabla \phi \cdot \nabla \hat{\phi} + \bar{\mu} \hat{\phi} \Big) \; dV, \quad \forall \hat{\phi}, \\
    \label{eq:stat_2}
    0 &= \int_{V} \Big( \frac{k}{\tilde{D}} \underline{j} \cdot \underline{\hat{j}} + \bar{\mu} \nabla \cdot \underline{\hat{j}} \Big) \; dV, \quad \forall \underline{\hat{j}}, \\
    \label{eq:stat_3}
    0 &= \int_{V} \Big( \frac{\partial \phi}{\partial t} + \nabla \cdot \underline{j} \Big) \hat{\bar{\mu}} \; dV, \quad \forall \hat{\bar{\mu}},
\end{align}

that is, the affinity relation \eqref{eq:Affinityvariational}, Onsager's law \eqref{eq:derived_j_} and the continuity equation \eqref{eq:continuity} in weak form. Only \eqref{eq:stat_3} contains the time derivative, and no time derivative of the stationarity conditions is taken (rule (b) of Section~\ref{sec:tep}). The time discretization is imposed after the variation, and only where a time derivative remains: in \eqref{eq:stat_3} the rate is replaced by the backward difference quotient $\partial \phi / \partial t \approx (\phi - \phi_t)/\Delta t$, with the time increment $\Delta t$, the previous solution $\phi_t = \phi(t)$ and the new solution $\phi = \phi(t + \Delta t)$, and all remaining quantities, in particular the driving force $\partial f / \partial \phi$ and the gradient $\nabla \phi$ in \eqref{eq:stat_1}, are evaluated at the new state (backward Euler). Thereby, time-dependent simulations in an implicit scheme are straightforward in their numerical implementation. The variational formulation for Galerkin-type methods follows by restricting the fields and the directions to finite element spaces:

\begin{align}
    \label{eq:var_form_1}
    0 &= \int_{V} \Big( \frac{\partial f(\phi)}{\partial \phi} \hat{\phi} + \kappa \nabla \phi \cdot \nabla \hat{\phi} + \bar{\mu} \hat{\phi} \Big) \; dV, \quad \forall \hat{\phi} \in W \\
    \label{eq:var_form_2}
    0 &= \int_{V} \Big( \frac{k}{\tilde{D}} \underline{j} \cdot \underline{\hat{j}} + \bar{\mu} \nabla \cdot \underline{\hat{j}} \Big) \; dV, \quad \forall \underline{\hat{j}} \in \mathbf{W}_0 \\
    \label{eq:var_form_3}
    0 &= \int_{V} \Big( \frac{\phi - \phi_t}{\Delta t} + \nabla \cdot \underline{j} \Big) \hat{\bar{\mu}} \; dV \quad \forall \hat{\bar{\mu}} \in W
\end{align}

Here $W$ and $\mathbf{W}$ denote the finite element spaces of the scalar and the vector-valued unknowns, specified below, and $\mathbf{W}_0$ the subspace of $\mathbf{W}$ whose elements satisfy the homogeneous boundary condition of the flux, $\underline{j} \cdot \underline{n} = 0$ on $\partial V$ for a closed system (in the computations the full flux vector is set to zero on $\partial V$, which implies it). Newer finite element tools like FEniCSx or COMSOL include automated variational differentiation of functionals, which allows for rapid implementation of models and simpler coupling of many phenomena, since only the Rayleighian has to be defined and the variational forms are computed automatically. This is especially useful for multi-component systems.

\textbf{Relation to the incremental minimization.} The discrete system \eqref{eq:var_form_1}--\eqref{eq:var_form_3} coincides with the optimality system of the constrained minimization of $\mathcal{F}[\phi] + \Delta t \cdot \tfrac{1}{2}\mathcal{P}[\underline{j}]$ over all pairs $(\phi, \underline{j})$ that satisfy the discrete mass balance $(\phi - \phi_t)/\Delta t + \nabla \cdot \underline{j} = 0$, with $\bar{\mu}$ as the multiplier of the constraint. Its Lagrangian, divided by $\Delta t$ and shifted by the constant $\mathcal{F}[\phi_t]$, is the time-discrete functional

\begin{equation}
    \label{eq:incremental}
    \mathcal{I}_{\Delta t}[\phi, \underline{j}, \bar{\mu}] := \frac{\mathcal{F}[\phi]-\mathcal{F}[\phi_t]}{\Delta t} + \int_V \Big( \frac{k}{2\tilde{D}} |\underline{j}|^2 +  \bar{\mu} \big( \frac{\phi - \phi_t}{\Delta t} + \nabla \cdot \underline{j}\big) \Big) \; dV,
\end{equation}

whose derivative with respect to $\phi$ in the direction $\hat{\phi}$ is $\frac{1}{\Delta t}\int_V (\partial_\phi f(\phi)\, \hat{\phi} + \kappa \nabla \phi \cdot \nabla \hat{\phi} + \bar{\mu} \hat{\phi})\, dV$, that is, \eqref{eq:var_form_1} up to the factor $1/\Delta t$, while its derivatives with respect to $\underline{j}$ and $\bar{\mu}$ are \eqref{eq:var_form_2} and \eqref{eq:var_form_3}. This is not a second derivation, since the Rayleighian is varied with respect to the rate and not with respect to the state, but a property of the discrete equations with two consequences. First, the backward Euler step is the minimizing movement, or incremental minimization, of the free energy in the metric defined by the dissipation \cite{jordan1998variational, otto2001geometry, flachberger2026variational}, the discrete counterpart of the statement that the physical rates minimize the Rayleighian at every instant. Second, since the objective is convex for a convex $f$ and coercive in $\underline{j}$, and the pair $(\phi_t, \underline{0})$ is admissible, the discrete problem has a solution for every $\Delta t > 0$, unique in $(\phi, \underline{j})$. The published code exploits the same fact: the notebook that produced the figures of this work defines the functional \eqref{eq:incremental} and lets FEniCSx derive its optimality conditions, which are exactly \eqref{eq:var_form_1}--\eqref{eq:var_form_3} (Section~\ref{sec:results}).

\textbf{Discrete thermodynamic consistency.} The structure of the scheme implies the three properties that Section~\ref{sec:degeneracy} asked for, independently of the choice of the finite element spaces and without any requirement on the regularity of $\phi$ beyond square integrability for $\kappa = 0$ (and $H^1$ for $\kappa > 0$).

\begin{proposition}
\label{prop:consistency}
Let $W$ contain the constant functions, let $\mathbf{W} \subset H(\mathrm{div}, V)$, let $\phi^n \in W$, let $f$ be convex (cases (i) and (ii), with $\kappa \geq 0$), and let $(\phi^{n+1}, \underline{j}^{n+1}, \bar{\mu}^{n+1})$ solve \eqref{eq:var_form_1}--\eqref{eq:var_form_3} with $\phi_t = \phi^n$. Then
\begin{enumerate}
    \item[(a)] $\int_V \phi^{n+1} \, dV = \int_V \phi^n \, dV$ (exact mass conservation);
    \item[(b)] $\mathcal{F}[\phi^{n+1}] + \Delta t \, \mathcal{P}[\underline{j}^{n+1}] + \mathcal{B}_n = \mathcal{F}[\phi^n]$, where
    \begin{equation*}
        \mathcal{B}_n := \int_V \Big( f(\phi^n) - f(\phi^{n+1}) - \frac{\partial f(\phi^{n+1})}{\partial \phi} \, (\phi^n - \phi^{n+1}) \Big) \, dV + \frac{\kappa}{2} \int_V |\nabla (\phi^{n+1} - \phi^n)|^2 \, dV \; \geq \; 0,
    \end{equation*}
    so that $\mathcal{F}[\phi^{n+1}] - \mathcal{F}[\phi^n] \leq - \Delta t \, \mathcal{P}[\underline{j}^{n+1}] \leq 0$ (discrete second law).
\end{enumerate}
\end{proposition}

\begin{proof}
(a) follows from \eqref{eq:var_form_3} with $\hat{\bar{\mu}} = 1$ and the divergence theorem, since $\underline{j}^{n+1} \cdot \underline{n} = 0$ on $\partial V$. For (b), test \eqref{eq:var_form_1} with $\hat{\phi} = \phi^{n+1} - \phi^n$, \eqref{eq:var_form_3} with $\hat{\bar{\mu}} = \bar{\mu}^{n+1}$ and \eqref{eq:var_form_2} with $\underline{\hat{j}} = \underline{j}^{n+1}$, which is admissible because $\underline{j}^{n+1} \in \mathbf{W}_0$:
\begin{align*}
    \int_V \Big( \frac{\partial f(\phi^{n+1})}{\partial \phi} (\phi^{n+1} - \phi^n) + \kappa \nabla \phi^{n+1} \cdot \nabla (\phi^{n+1} - \phi^n) \Big) dV &= - \int_V \bar{\mu}^{n+1} (\phi^{n+1} - \phi^n) \, dV \\
    &= \Delta t \int_V \bar{\mu}^{n+1} \nabla \cdot \underline{j}^{n+1} \, dV \\
    &= - \Delta t \int_V \frac{k}{\tilde{D}} |\underline{j}^{n+1}|^2 \, dV .
\end{align*}
The right-hand side is $-\Delta t \, \mathcal{P}[\underline{j}^{n+1}]$. On the left-hand side, the identity $\kappa \nabla a \cdot \nabla (a - b) = \tfrac{\kappa}{2} \big( |\nabla a|^2 - |\nabla b|^2 + |\nabla (a - b)|^2 \big)$ and the definition of $\mathcal{B}_n$ give $\mathcal{F}[\phi^{n+1}] - \mathcal{F}[\phi^n] + \mathcal{B}_n$. The first integral of $\mathcal{B}_n$ is the Bregman divergence of the convex function $f$ and is nonnegative, as is the second.
\end{proof}

Three remarks are in order. First, $\mathcal{B}_n$ is the numerical dissipation of the backward Euler step; it is of second order in the increment $\phi^{n+1} - \phi^n$ and vanishes as $\Delta t \to 0$, and the discrete energy identity closes exactly when it is accounted for, as in the computations of \cite{flachberger2026variational}. Second, the proof uses only that $\phi^{n+1} - \phi^n$, $\underline{j}^{n+1}$, $\bar{\mu}^{n+1}$ and the constants are admissible test functions; it holds for every choice of conforming finite element spaces, continuous or discontinuous, and for every convex $f$, so that $f'' = 0$ on an interval, the degenerate case (ii), is explicitly allowed. The gradient term may moreover be replaced by any symmetric positive semidefinite quadratic form $a(\phi, \phi)$, such as the facet jump penalty of Section~\ref{sec:spaces} for discontinuous spaces, with $\tfrac{1}{2} a(\phi^{n+1} - \phi^n, \phi^{n+1} - \phi^n)$ in place of the gradient term of $\mathcal{B}_n$. This is the precise sense in which the scheme operates on the boundary of degeneracy: mass conservation, energy bookkeeping and the second law do not depend on the parabolicity of the equation. Third, for the double-well $f_0$ the Bregman term is of indefinite sign in the spinodal region, and the inequality then holds up to a term quadratic in the increment, which is the general situation for the fully implicit Cahn--Hilliard scheme; convex splitting of $f_0$ would restore it unconditionally.

\textbf{Choice of the finite element spaces.} Forms \eqref{eq:var_form_2} and \eqref{eq:var_form_3} are the dual mixed formulation of the second-order problem \eqref{eq:cahn-hilliard}, introduced for elliptic problems by Raviart and Thomas \cite{raviart1977mixed} and Brezzi, Douglas and Marini \cite{brezzi1985two}, see \cite{boffi2013mixed}: by solving for the flux explicitly, the second-order elliptic operator is split into two first-order equations, and the derivative is thrown onto the flux test function. Mixed methods of this type have been used for degenerate parabolic equations arising in porous media flow \cite{arbogast1996nonlinear}. The innovative aspect of the present formulation lies in its treatment of the affinity: no differential operator acts on $\bar{\mu}$, and, for $\kappa = 0$, none acts on $\phi$, in contrast to formulations in which the mole fraction is the only degree of freedom \cite{svoboda2013new, svoboda2017incorporation, flachberger2024numerical}. The natural function spaces of the dual mixed form are $\underline{j} \in H(\mathrm{div}, V)$ and $\phi, \bar{\mu} \in L^2(V)$ (for $\kappa = 0$), so that a discontinuous mole fraction is admissible and the convex free energy of Figure~\ref{fig1} can be used without regularization, i.e.\ case (ii) can be implemented. In \cite{flachberger2024numerical}, a detailed analysis reveals why traditional variational forms for reactive diffusion problems fail due to the continuity properties of the affinity; by avoiding differentiation of the affinity, our approach circumvents these issues and is able to solve the cases discussed.

In the computations reported below we use, for all three fields, equal-order continuous piecewise quadratic elements on a triangulation $\mathcal{T}_h$ of $V$ (intervals in one dimension, quadrilaterals in two), $W = H^1(V) \cap P_2(\mathcal{T}_h)$ for $\hat{\phi}$ and $\hat{\bar{\mu}}$ and $\mathbf{W} = [H^1(V)]^d \cap [P_2(\mathcal{T}_h)]^d$ for $\underline{\hat{j}}$ ($Q_2$ on quadrilaterals), with the full flux vector prescribed to vanish on $\partial V$, which implies $\underline{j} \cdot \underline{n} = 0$. This is a pragmatic choice in view of the mechanical coupling we aim at, where the displacement is discretized with continuous elements, and it allows the scheme to be carried over to established finite element frameworks such as ABAQUS, which provide continuous elements only. Two facts should be stated. For the stationary Darcy-type subproblem \eqref{eq:var_form_2}--\eqref{eq:var_form_3} alone, equal-order continuous pairs do not satisfy the inf-sup condition \cite{masud2002stabilized}, whereas the Raviart--Thomas or Brezzi--Douglas--Marini spaces for the flux together with discontinuous spaces for the scalars do \cite{boffi2013mixed}. In the time-dependent problem, the term $(\phi - \phi_t)/\Delta t$ of \eqref{eq:var_form_3} couples the multiplier to the mole fraction through \eqref{eq:var_form_1} and contributes, in the bulk phases, a term in $\bar{\mu}$ of the stabilizing sign, proportional to $1/(f'' \Delta t)$, that the stationary problem lacks; it is this term that keeps the equal-order discretization usable at finite $\Delta t$. In the two-phase region, where $f'' = 0$, equation \eqref{eq:var_form_1} no longer involves $\phi$ at all: the affinity is pinned to its common-tangent value, and $\phi$ is determined by the mass balance \eqref{eq:var_form_3} alone. This is the indeterminacy of the degenerate region, not a failure of the inf-sup condition, and it is the origin of the fluctuations discussed in Section~\ref{sec:results}. We make no claim of uniform inf-sup stability for the equal-order choice; Proposition~\ref{prop:consistency} does not depend on it, and the inf-sup stable Raviart--Thomas pair is available as an alternative within the same formulation (Section~\ref{sec:spaces}).

\subsection{A ternary system with sources and sinks for vacancies}
\label{subsec:ternary}

Finally, we want to give an example of how to implement a specific model
within the presented methodology by considering the vacancy diffusion
model derived by Svoboda in \cite{svoboda2006diffusion}. The concept of
``vacancy'' refers to empty or vacant sites in a crystal lattice, which
are crucial for enabling diffusion in alloys. By allowing atoms to
switch positions with vacant sites, vacancies facilitate the movement of
atoms within an alloy. Since vacancies are merely empty sites in a
crystal lattice, it is essential to recognize that they are purely
geometric features and do not constitute matter. This distinction is
important, as it permits us to treat vacancies as a non-conserved
component, allowing for more nuanced modeling of diffusion processes.
Therefore Svoboda suggested treating a binary system with vacancies like
a ternary system (with vacancies constituting the third component). The
composition of the system is then defined by so-called site fractions
instead of mole fractions, which refer to the ratio of the amount of a
component to the number of available lattice sites in the
atomic grid. We begin by defining the free energy of the system
which will just be a superposition of the convex hull from Figure~\ref{fig1},
namely \(f_{0}^{**}\) and a simple quadratic contribution of the site
fraction of vacancies \(\phi_0\):

\begin{equation}
    \label{eq:free}
    \mathcal{F}[\phi_0, \phi_1] := \int_V f_m \; dV = \int_V \Big( f_{0}^{**}(\phi_1) + \tfrac{\kappa}{2} | \nabla \phi_1 |^2 + \tfrac{k_0}{2} \cdot (\phi_0 - \phi_{0}^{eq})^2 \Big) \; dV
\end{equation}

As can be seen, the free energy will consist only of the molar free
energy of the equilibrium system \(f_{0}^{**}\) if the vacancies \(\phi_0\)
assume their equilibrium site fraction \(\phi_{0}^{eq}\). A perfect crystal
without vacancies and dislocations may possess less internal energy than
one with imperfections; however, the presence of imperfections increases
the configurational entropy of the crystal, thereby reducing the free
energy. The equilibrium site fraction of vacancies marks the minimum of
the overall free energy of the system with respect to vacancies. The dissipation, one half of which enters the Rayleighian \eqref{eq:rayleighian_def}, is defined according to \cite{svoboda2006diffusion} by the following functional:

\begin{equation}
    \mathcal{P}[\underline{j}_0, \underline{j}_1, \psi] := \int_V \mathscr{P} \; dV = \int_V \Big( \frac{k_0}{D_0} |\underline{j}_0|^2 + \frac{k_1}{D_1} |\underline{j}_1|^2 + \frac{k_2}{D_2} |-\underline{j}_0-\underline{j}_1|^2 + \frac{1}{A_{\psi}} \psi^2 \Big) \; dV
\end{equation}

It is worth noting that the diffusion coefficients of individual components,
$D_i$, are utilized in this formulation. Furthermore, through the definition of the flux of component 2, given by equation \eqref{eq:flux_constraint_2} as $j_2:=-j_0-j_1$, the fluxes become intertwined according to Onsager's theory. This coupling not only constrains the system and enforces the vacancy mechanism but also gives rise to non-diagonal entries in the resulting Onsager coefficient matrix of the system.
The Onsager coefficients can be calculated as follows (in the scaled form of Section~\ref{sec:scaling}; in the unscaled form a factor $T$ appears):

\begin{equation}
    \label{eq:onsgeronsger}
    \tilde{L}_{ik} = \Big[ \Big( \frac{\partial^2 \tfrac{1}{2}\mathscr{P}}{\partial \underline{j}_i \partial \underline{j}_k} \Big)_{T,V,\{\phi\}} \Big]^{-1},  \quad  \forall i,k\in\{1,\ldots,n-1\}
\end{equation}

Consequently, the fluxes are coupled, and a driving force for one component can induce fluxes of other components, as predicted by Onsager's reciprocal relations.
The function \(\psi\) is characteristic of this
 model and represents the rate of generation and annihilation of
vacancies. Since vacancies are a non-conserved component they can emerge
and vanish in the bulk to facilitate diffusion. Their evolution in
time is therefore also not determined by a pure conservation equation but
rather a transport-type equation. The evolution equations are considered
in the constraint functional $\mathcal{C} = \int_V \mathscr{C} \; dV$, whose density in the rate form of Section~\ref{sec:tep} reads, with the coefficients of $\psi$ evaluated at the given state (rule (c)),

\begin{multline}
    \label{eq:constraint_ternary}
    \mathscr{C}\Big(\frac{\partial \phi_0}{\partial t}, \frac{\partial \phi_1}{\partial t}, \underline{j}_0, \underline{j}_1, \bar{\mu}_0, \bar{\mu}_1, \psi\Big) := \\
    \bar{\mu}_0 \Big( \frac{\partial \phi_0}{\partial t} + \nabla \cdot \underline{j}_0 - \psi \cdot (1-\phi_{0}) \Big) + \bar{\mu}_1 \Big( \frac{\partial \phi_1}{\partial t} + \nabla \cdot \underline{j}_1 + \psi \cdot \phi_{1} \Big) .
\end{multline}

As can be seen, the rate of annihilation and generation of vacancies
\(\psi\) not only impacts the evolution equation of vacancies but also
of component \(\phi_1\). This is because an increase of vacancies also
increases the number of available sites in the lattice which impacts all
site fractions. Therefore, the generation of vacancies also causes a
strain and a change in the molar volume \(\Omega\). Strictly, a lattice whose number of sites changes is no longer at rest, and rule (a) of Section~\ref{sec:tep} would call for the material form of the transport theorem together with the mechanical treatment of the problem. This was
nevertheless neglected since it would also require the mechanical
treatment of the problem, which is beyond the scope of this publication.
The final problem is defined by its Rayleighian density, the rate of the free energy \eqref{eq:free} in the direction of the rates, plus one half of the dissipation, plus the constraint:

\begin{multline}
    \label{eq:rayleighian_ternary_density}
    \mathscr{R}\Big(\frac{\partial \phi_0}{\partial t}, \frac{\partial \phi_1}{\partial t}, \underline{j}_0, \underline{j}_1, \bar{\mu}_0, \bar{\mu}_1, \psi\Big) := \frac{\partial f_m}{\partial \phi_0}\, \frac{\partial \phi_0}{\partial t} + \frac{\delta f_m}{\delta \phi_1}\, \frac{\partial \phi_1}{\partial t} + \tfrac{1}{2}\mathscr{P}(\underline{j}_0, \underline{j}_1, \psi) \\
    + \mathscr{C}\Big(\frac{\partial \phi_0}{\partial t}, \frac{\partial \phi_1}{\partial t}, \underline{j}_0, \underline{j}_1, \bar{\mu}_0, \bar{\mu}_1, \psi\Big)
\end{multline}

The variational form is determined by computing the first variation with respect to the rates $\partial \phi_0 / \partial t$, $\partial \phi_1 / \partial t$, the fluxes $\underline{j}_0, \underline{j}_1$, the multipliers $\bar{\mu}_0, \bar{\mu}_1$ and the rate $\psi$, at the given state $(\phi_0, \phi_1)$:

\begin{align}
    \label{eq:lagrangian_ternary}
    \delta \int_{V} \mathscr{R}\Big(\frac{\partial \phi_0}{\partial t}, \frac{\partial \phi_1}{\partial t}, \underline{j}_0, \underline{j}_1, \bar{\mu}_0, \bar{\mu}_1, \psi\Big) \; dV = 0 .
\end{align}

After the variation, the time discretization of Section~\ref{subsec:fem} is imposed: the rates are replaced by the backward difference quotients $(\phi_0 - \phi_0^t)/\Delta t$ and $(\phi_1 - \phi_1^t)/\Delta t$, the driving forces are evaluated at the new state, and the coefficients $1 - \phi_0$ and $\phi_1$ of $\psi$ are taken at the previous time level, so that the discrete constraint equations read

\begin{equation}
    \label{eq:constraint_ternary_discrete}
    \int_V \Big( \frac{\phi_0 - \phi_{0}^{t}}{\Delta t} + \nabla \cdot \underline{j}_0 - \psi \, (1-\phi_{0}^{t}) \Big) \hat{\bar{\mu}}_0 \; dV = 0, \qquad
    \int_V \Big( \frac{\phi_1 - \phi_{1}^{t}}{\Delta t} + \nabla \cdot \underline{j}_1 + \psi \, \phi_{1}^{t} \Big) \hat{\bar{\mu}}_1 \; dV = 0 \qquad \forall \hat{\bar{\mu}}_0, \hat{\bar{\mu}}_1 .
\end{equation}

\subsection{Variable diffusion coefficients}

It is commonly observed that the components of an alloy exhibit different diffusion coefficients in different phases. The diffusion coefficients $D_i$, or more generally the Onsager coefficients $\tilde{L}_{ik}$, may therefore be chosen phase-dependent, that is, as functions of the mole fraction $\phi$. This requires no special treatment. The Rayleighian is varied with respect to the rates, the fluxes and the multipliers at the given state (rule (c) of Section~\ref{sec:tep}), so that $\phi$ is not touched by the variation, and coefficients depending on it enter the stationarity conditions, and the discrete forms \eqref{eq:var_form_1}--\eqref{eq:var_form_3}, simply as given functions of the state. In the implicit scheme they may be evaluated at the new state or, for simplicity, at the previous time level.

\section{Results}
\label{sec:results}

To explore the capabilities of our method, we employed the open-source
finite element solver FEniCSx, which is continuously developed by the
FEniCS-community. The software relies on several key modules, including
dolfinx \cite{barrata2023dolfinx}, basix
\cite{scroggs2022construction, scroggs2022basix}, and ufl \cite{alnaes2015fenics}.
FEniCSx
was used in conjunction with Matplotlib \cite{hunter2007matplotlib} to generate all graphics.
The codes to run the simulations are publicly available\footnote{\url{https://github.com/wolfgang-alpha/Modelling-1st-and-2nd-Order-Phase-Transformations}} and the automated functional derivative feature of FEniCSx is also utilized to generate the variational forms for the ternary system defined by equations \eqref{eq:free}-\eqref{eq:lagrangian_ternary}. The published code provides both constructions of Section~\ref{subsec:fem}: the notebook that produced the figures of this work defines the incremental functional \eqref{eq:incremental} and lets FEniCSx derive its stationarity conditions, the Raviart--Thomas notebook of Section~\ref{sec:spaces} generates the discrete forms by variation with respect to the rates and multipliers at the frozen state followed by the backward Euler substitution, and a script verifies that the two constructions yield identical forms, to machine precision, for both element choices. A second script checks the mass balance and the identity of Proposition~\ref{prop:consistency}(b) numerically, see below. All parameters are given in the scaled units of Section~\ref{sec:scaling}, with lengths in mm, times in s and molar free energies in J/mol; the gradient coefficient $\kappa$ is quoted as it enters $f_m$ in \eqref{eq:case_ch}.

To showcase the capabilities of our method, we start by investigating one- and two-dimensional
variational problems defined by forms \eqref{eq:var_form_1} - \eqref{eq:var_form_3}. First
the one-dimensional case is explored, where we contrast solutions for the classical
Cahn-Hilliard model and the cases (i) and (ii), all initialized
with a common condition, namely a linearly rising distribution of the
mole fraction.

\begin{figure}[h!]
    \centering
    \includegraphics[max size={0.7\linewidth}{0.7\paperheight}]{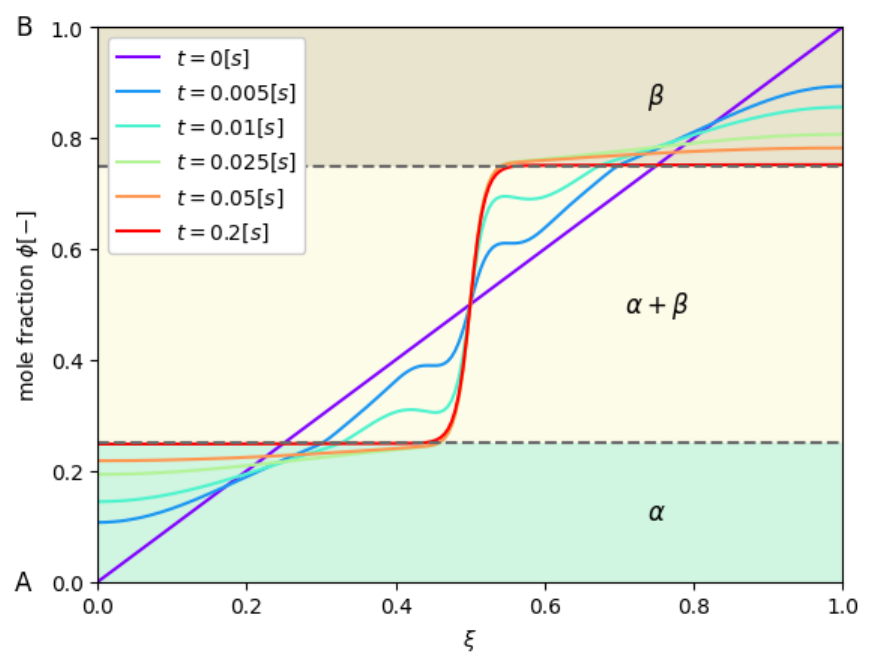}
    \caption{System evolution of the variational forms given in \eqref{eq:var_form_1} - \eqref{eq:var_form_3}.
The contour visible is the distribution of the mole fraction \(\phi\) after
different times. The molar free energy \(f_0(\phi)\) from Figure~\ref{fig1} is
employed (analytical expression for $f_0(\phi)$ can be found in Appendix~\ref{app:functions}). The homogenization parameter was set to \(\kappa = 10^{-4} [\frac{J mm^2}{mol}]\), the inter-diffusion coefficient to \(\tilde{D} = 2[mm^2/s]\) and $k=2[\tfrac{J}{mol}]$.}
    \label{fig2}
    \end{figure}

\begin{figure}[H]
    \centering
    \includegraphics[max size={0.7\linewidth}{0.7\paperheight}]{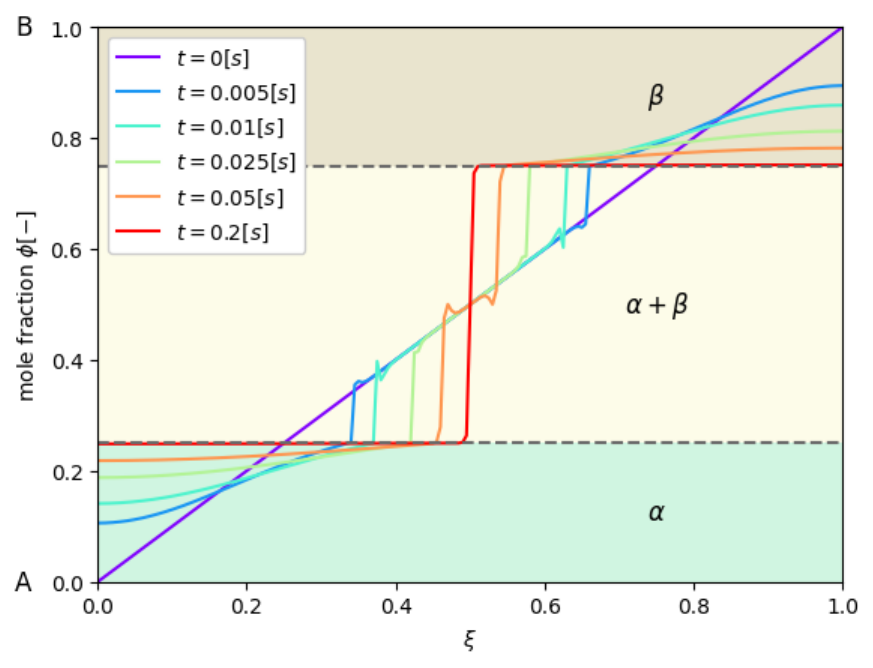}
    \caption{System evolution of the variational forms given in \eqref{eq:var_form_1} - \eqref{eq:var_form_3} for case (ii). The contour
visible is the distribution of the mole fraction \(\phi\) after different
times. The convex molar free energy \(f_{0}^{**}(\phi)\) from Figure~\ref{fig1} is
employed (analytical expression for $f_{0}^{**}(\phi)$ can be found in Appendix~\ref{app:functions}). The homogenization parameter was set to $\kappa = 0$, the inter-diffusion coefficient to \(\tilde{D} = 2[mm^2/s]\) and $k=2[\tfrac{J}{mol}]$. The fluctuations at the moving edges of the two-phase region are the signature of the degeneracy discussed in Section~\ref{sec:degeneracy}; they vanish together with the region.}
    \label{fig3}
    \end{figure}

\begin{figure}[H]
    \centering
    \includegraphics[max size={0.7\linewidth}{0.7\paperheight}]{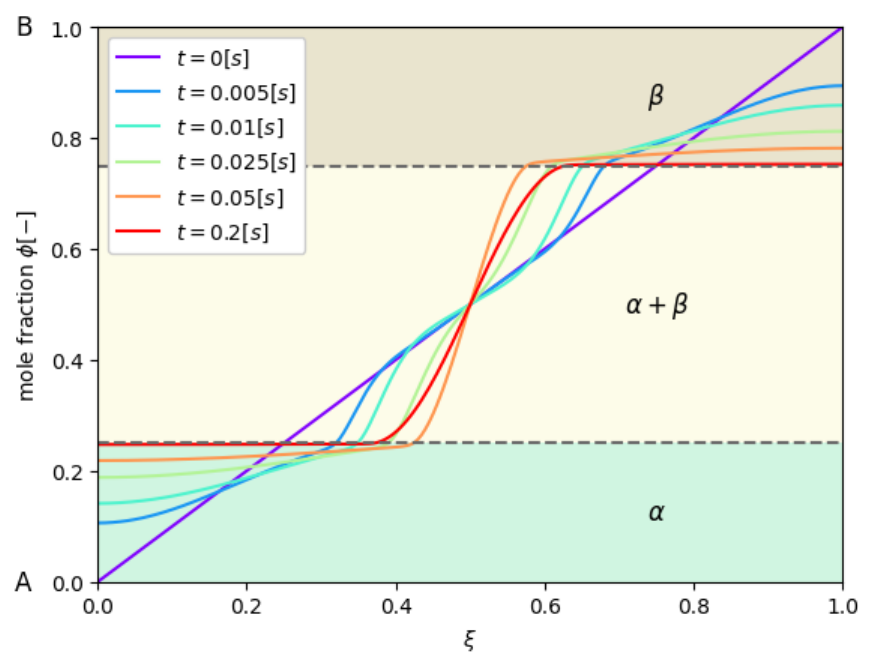}
    \caption{System evolution of the variational forms given in \eqref{eq:var_form_1} - \eqref{eq:var_form_3} for case (i). The contour
visible is the distribution of the mole fraction \(\phi\) after different
times. The convex molar free energy \(f_{0}^{**}(\phi)\) from Figure~\ref{fig1} is
employed (analytical expression for $f_{0}^{**}(\phi)$ can be found in Appendix~\ref{app:functions}). The homogenization parameter was set to \(\kappa = 10^{-4} [\frac{J mm^2}{mol}]\), the inter-diffusion coefficient to \(\tilde{D} = 2[mm^2/s]\) and $k=2[\tfrac{J}{mol}]$.}
    \label{fig4}
    \end{figure}

Figures~\ref{fig2} to~\ref{fig4} illustrate the time dependent phase
distributions obtained for each scenario. As anticipated, the classical
Cahn-Hilliard model exhibits phase separation, characterized by a smooth
transition between the phases (Figure~\ref{fig2}).

In contrast, special case (ii)
features a sharp, discontinuous interface due to the absence of
homogenization terms and the presence of a convex hull curve, as shown
in Figure~\ref{fig3}. The initial condition lies, for $0.25 < \xi < 0.75$, entirely inside the two-phase interval of the convex hull: this is the degenerate region of Section~\ref{sec:degeneracy}, in which the affinity is constant and no driving force exists. The region is consumed from both sides, as the bulk phases draw matter across its boundaries, until only a sharp interface at $\xi = 0.5$ remains, the position dictated by mass conservation. The fluctuations visible at the moving edges of the two-phase region before the interface has sharpened are the signature of the degeneracy, not of an instability of the scheme: inside the region the equation provides no damping of perturbations of the mole fraction, since any profile with values in $[\phi_\alpha, \phi_\beta]$ is at rest there, and the discrete problem inherits this indeterminacy at the scale of the mesh. In all our computations the fluctuations remain inside the two-phase region, they disappear with it, and Proposition~\ref{prop:consistency} guarantees that mass is conserved and the free energy decreases monotonically throughout. We have verified this for the one-dimensional computations of Figures~\ref{fig2} to~\ref{fig4}: at every time step the identity of Proposition~\ref{prop:consistency}(b) closes to machine precision ($|\mathcal{F}[\phi^{n+1}] + \Delta t\,\mathcal{P}[\underline{j}^{n+1}] + \mathcal{B}_n - \mathcal{F}[\phi^n]| < 10^{-15}$ in the scaled units), the total mass drifts by less than $10^{-15}$, the numerical dissipation $\mathcal{B}_n$ is nonnegative, also for the double well, and the free energy decreases monotonically; the script performing this check is published with the code. The method operates, by design, on the boundary of degeneracy.
By incorporating an added
regularization, special case (i) yields a smooth phase transition,
similar to that observed for the classical Cahn-Hilliard model, but with
a wider phase transition region, as depicted in Figure~\ref{fig4}.

\begin{figure}[H]
    \centering
    \includegraphics[max size={0.7\linewidth}{0.7\paperheight}]{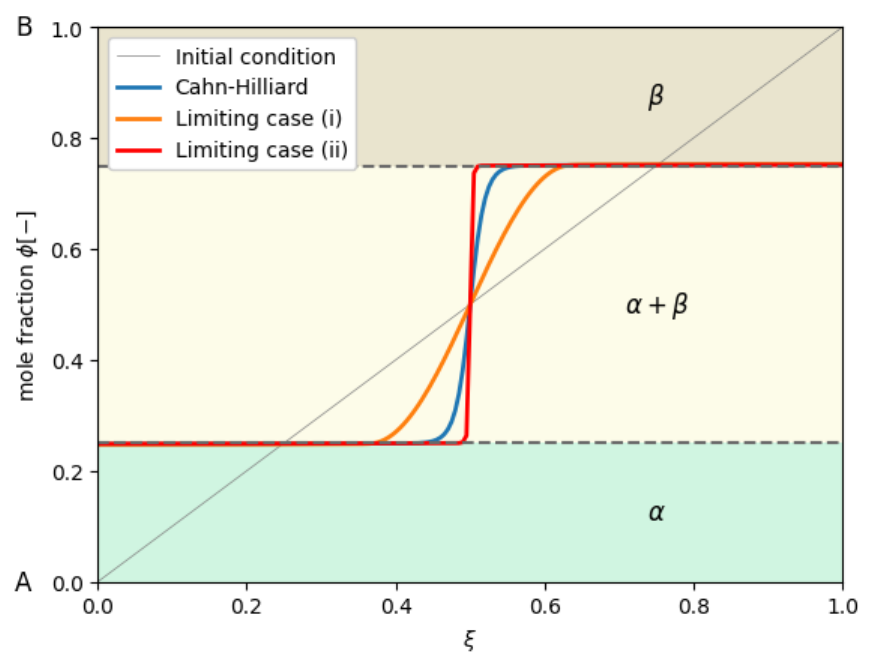}
    \caption{Direct comparison of the final states (\(t=0.2[s]\))
of the cases (i) and (ii) as well as the Cahn-Hilliard model, for the computations from Figures~\ref{fig2}, \ref{fig3} and \ref{fig4}.}
    \label{fig5}
    \end{figure}

\begin{figure}[H]
    \centering
    \includegraphics[max size={0.7\linewidth}{0.7\paperheight}]{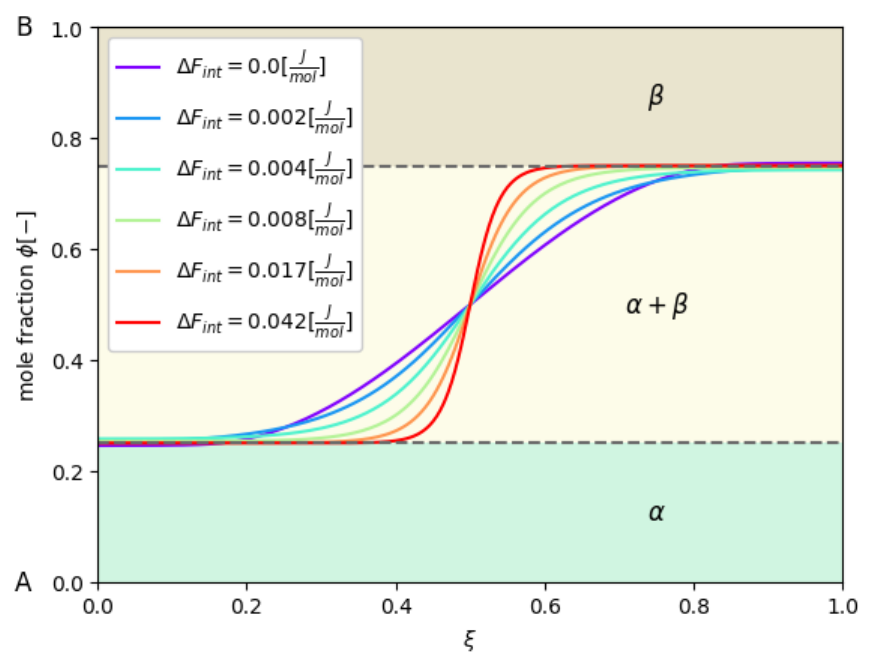}
    \caption{System evolution of the variational forms given in \eqref{eq:var_form_1} - \eqref{eq:var_form_3} for
different choices of the interface energy \(\Delta F_{int}\) (i.e.~the
relative height of the ``bump'' in the double-well). The contours
visible are the distributions of the mole fraction \(\phi\) after
\(0.2[s]\). The homogenization parameter was set to \(\kappa = 2 \cdot 10^{-3} [\frac{J mm^2}{mol}]\), the inter-diffusion coefficient to \(\tilde{D} = 2[mm^2/s]\) and $k=2[\tfrac{J}{mol}]$.}
    \label{fig6}
    \end{figure}

For comparison
the final states of the simulations are summarized in Figure~\ref{fig5}. The
disparity in the interface thickness of special case (i) and the
Cahn-Hilliard model can be attributed to the convex free energy not
imposing an energetic penalty on the interface. To confirm this, we have
performed the calculation again for many different interface energies.
In Figure~\ref{fig6} it can be seen that the interface energy clearly
impacts the width of the interface. In the classical Cahn-Hilliard
model, the ``bump'' in the double-well potential counteracts the
regularization term, resulting in a smooth yet locally restricted
interface.

\begin{figure}[H]
    \centering
    \includegraphics[max size={0.6\linewidth}{0.6\paperheight}]{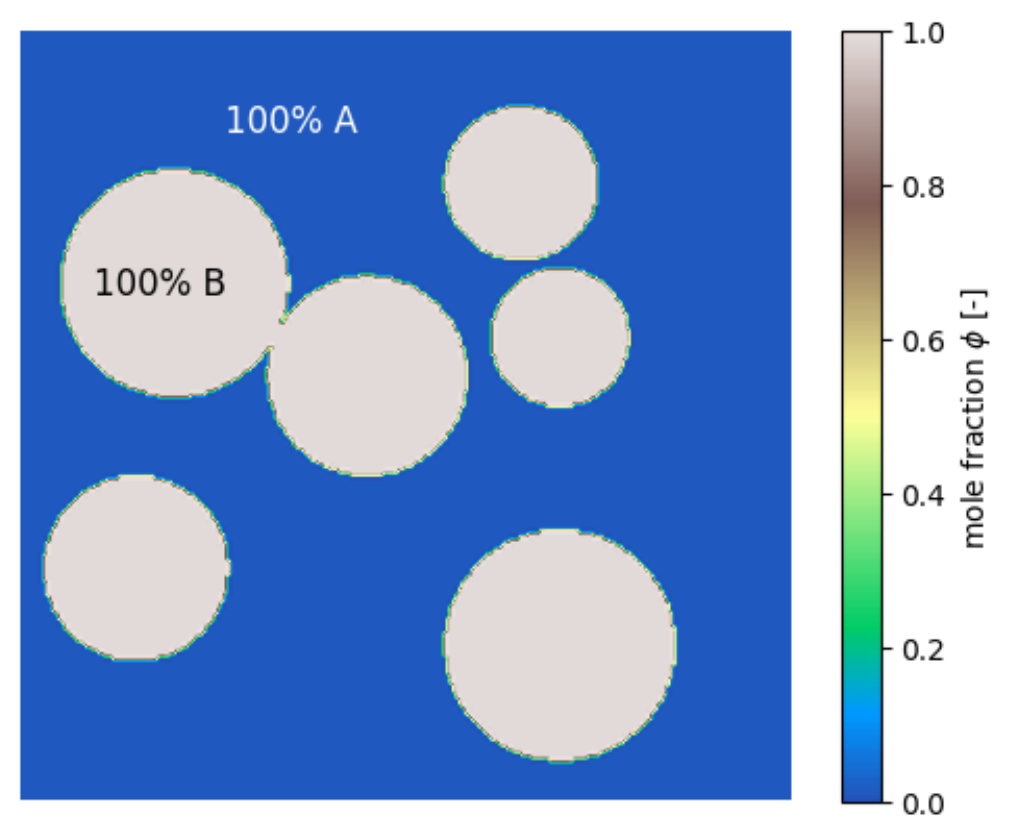}
    \caption{The initial condition of the mole fraction \(\phi\) in a
two-dimensional domain of unit length. The components of the binary
system are entirely separated from each other. The blue region
corresponds to component 2 and phase \(\alpha\) and the red region
corresponds to component 1 and phase \(\beta\).}
    \label{fig7}
    \end{figure}

\begin{figure}[H]
    \centering
    \includegraphics[max size={0.6\linewidth}{0.6\paperheight}]{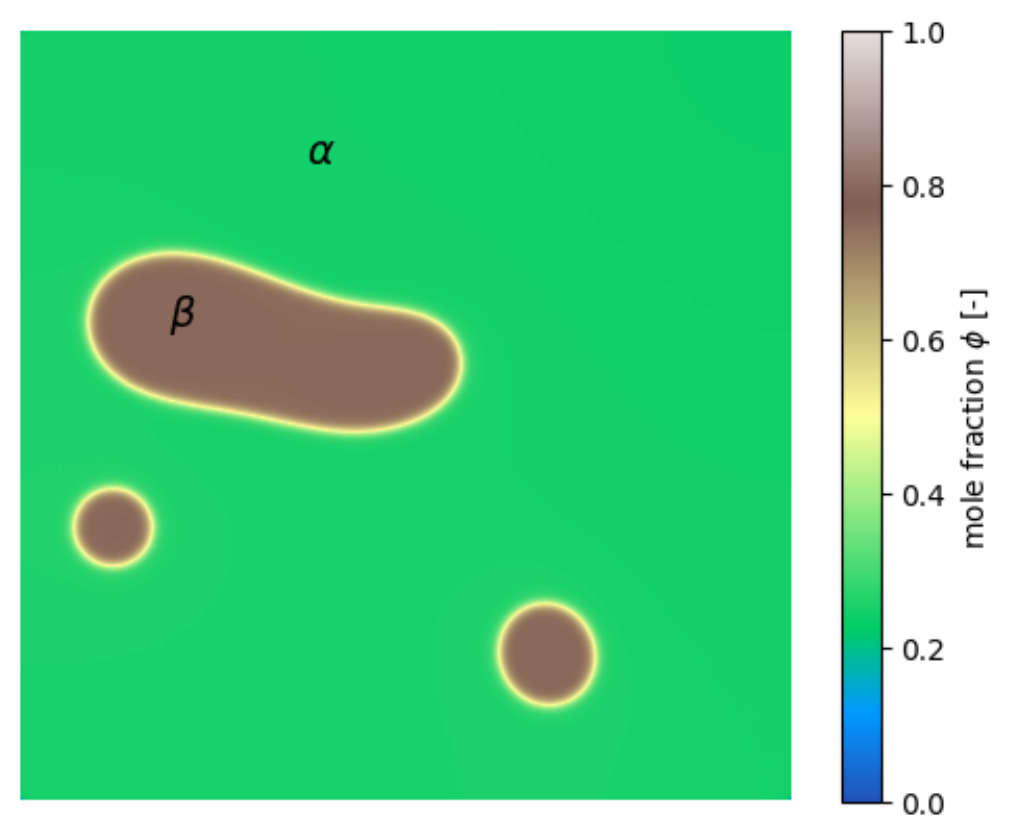}
    \caption{System evolution of the variational forms given in \eqref{eq:var_form_1} - \eqref{eq:var_form_3}.
The contour visible is the distribution of the mole fraction \(\phi\) after
\(0.1[s]\). The molar free energy \(f_0(\phi)\) from Figure~\ref{fig1} is employed
and the homogenization parameter was set to \(\kappa = 2 \cdot 10^{-5} [\frac{J mm^2}{mol}]\), the inter-diffusion coefficient to \(\tilde{D} = 2[mm^2/s]\) and $k=2[\tfrac{J}{mol}]$.}
    \label{fig8}
    \end{figure}

\begin{figure}[H]
    \centering
    \includegraphics[max size={0.6\linewidth}{0.6\paperheight}]{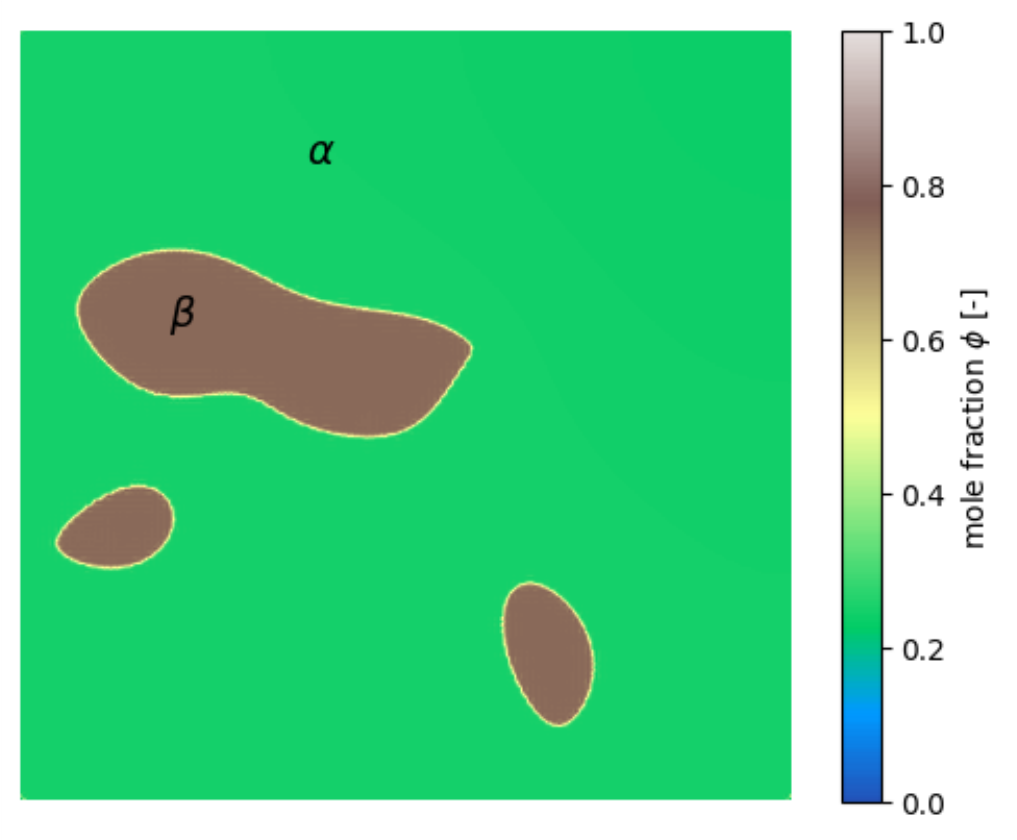}
    \caption{System evolution of the variational forms given in \eqref{eq:var_form_1} - \eqref{eq:var_form_3} for case (ii). The contour
visible is the distribution of the mole fraction \(\phi\) after \(0.1[s]\).
The convex molar free energy \(f_{0}^{**}(\phi)\) from Figure~\ref{fig1} is employed
and the homogenization parameter was set to \(\kappa = 0\), the inter-diffusion coefficient to \(\tilde{D} = 2[mm^2/s]\) and $k=2[\tfrac{J}{mol}]$.}
    \label{fig9}
    \end{figure}

\begin{figure}[H]
    \centering
    \includegraphics[max size={0.6\linewidth}{0.6\paperheight}]{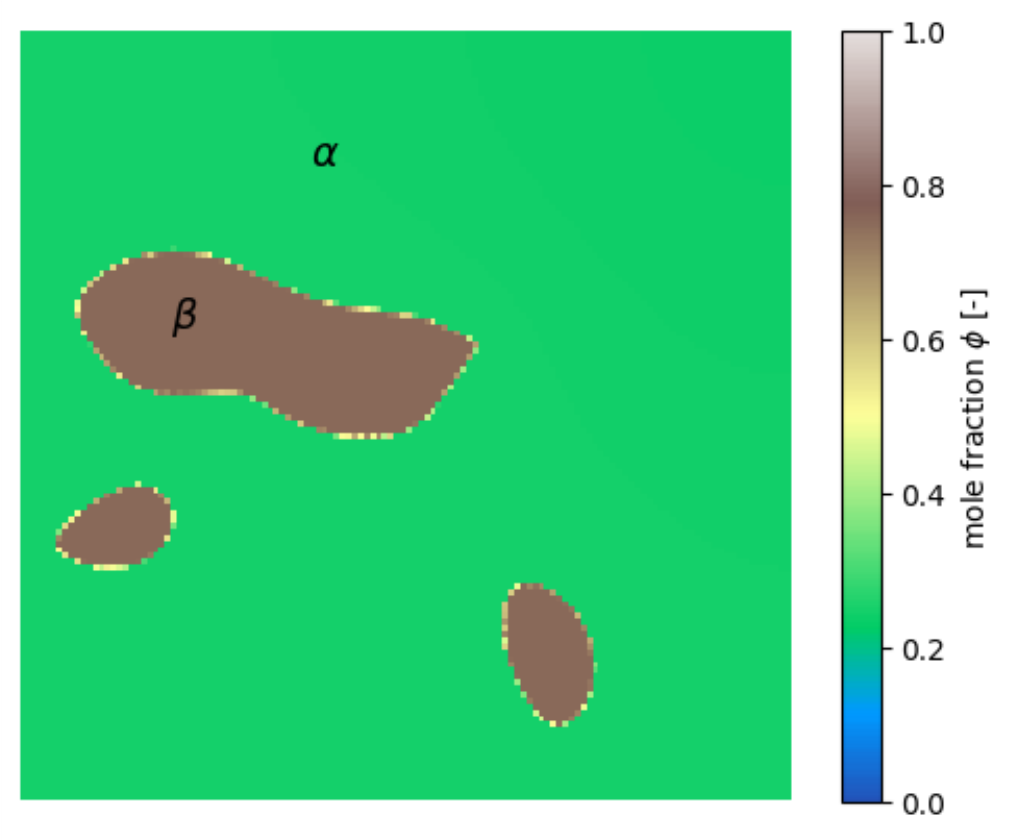}
    \caption{Cross-check of Figure~\ref{fig9} with the Raviart--Thomas/discontinuous pair of Section~\ref{sec:spaces} (RTCF$_1$ for the flux, piecewise constants for $\phi$ and $\bar{\mu}$) on the same mesh of $128 \times 128$ quadrilaterals, with the same parameters and time step; distribution of the mole fraction $\phi$ after $0.1[s]$. The morphology coincides with that of the equal-order computation of Figure~\ref{fig9}, and the interface is one cell wide. The computation took 4 minutes on a single thread, against 109 minutes for Figure~\ref{fig9}.}
    \label{fig9b}
    \end{figure}

\begin{figure}[H]
    \centering
    \includegraphics[max size={0.6\linewidth}{0.6\paperheight}]{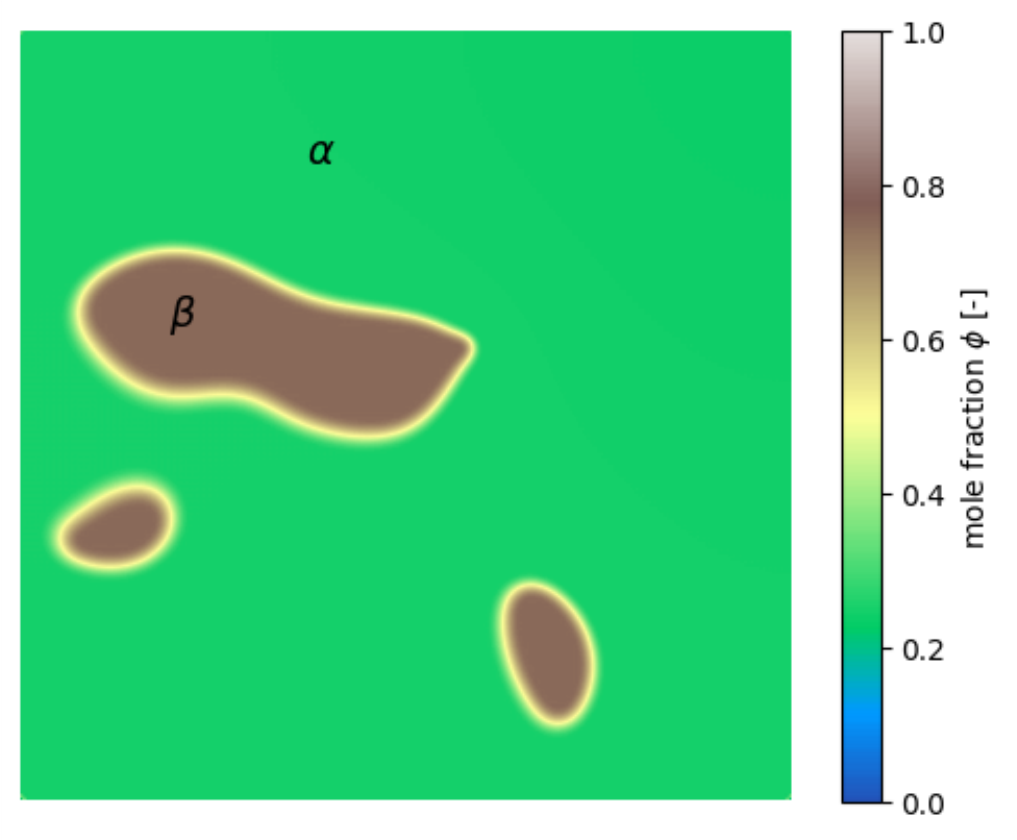}
    \caption{System evolution of the variational forms given in \eqref{eq:var_form_1} - \eqref{eq:var_form_3} for case (i). The contour
visible is the distribution of the mole fraction \(\phi\) after \(0.1[s]\).
The convex molar free energy \(f_{0}^{**}(\phi)\) from Figure~\ref{fig1} is employed
and the homogenization parameter was set to \(\kappa = 2 \cdot 10^{-7} [\frac{J mm^2}{mol}]\), the inter-diffusion coefficient to \(\tilde{D} = 2[mm^2/s]\) and $k=2[\tfrac{J}{mol}]$.}
    \label{fig10}
    \end{figure}

\begin{figure}[H]
    \centering
    \includegraphics[max size={0.6\linewidth}{0.6\paperheight}]{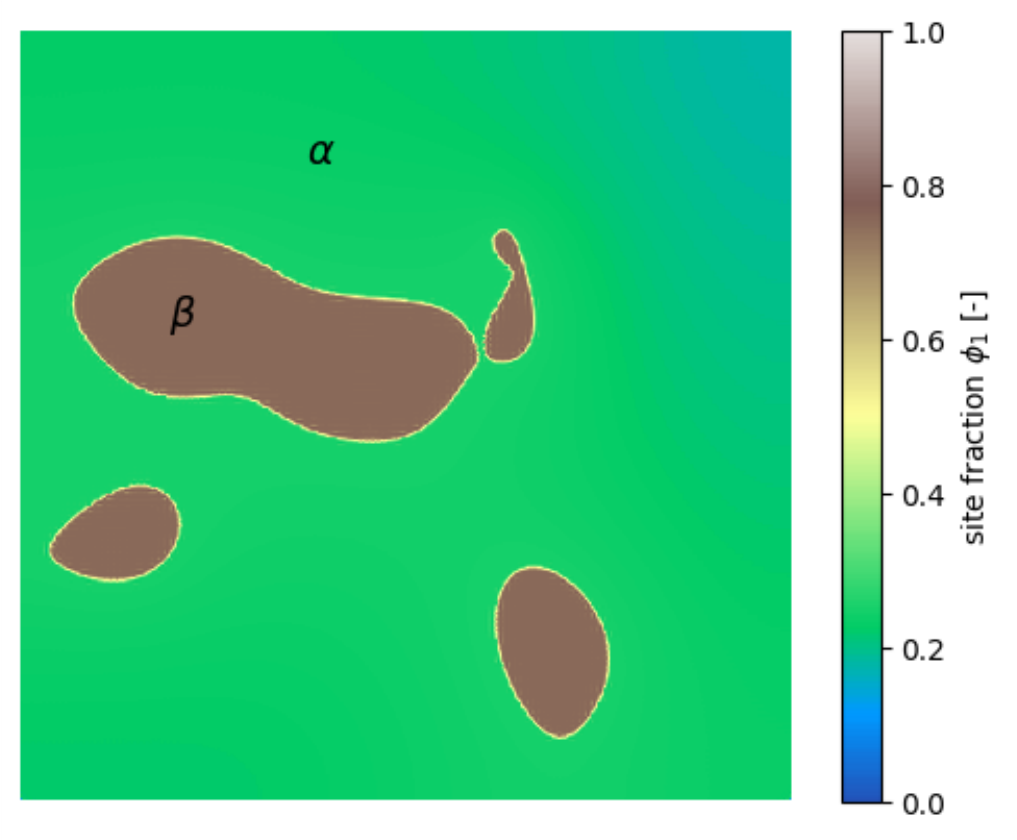}
    \caption{System evolution for the variational problem defined
in equation \eqref{eq:lagrangian_ternary}. The contour visible is the distribution of the site
fraction of component \(\phi_1\) after \(0.1[s]\). The parameters were
chosen as follows: \(\phi_{0}^{eq}=0.001\),
\(D_0 \rightarrow \infty\), \(D_1=2\), \(D_2=1\), $k_0=200[\tfrac{J}{mol}]$, $k_1 = k_2 = 2[\tfrac{J}{mol}]$, $A_\psi = 1/2$ and \(\kappa = 0\).}
    \label{fig11}
    \end{figure}

\begin{figure}[H]
    \centering
    \includegraphics[max size={0.6\linewidth}{0.6\paperheight}]{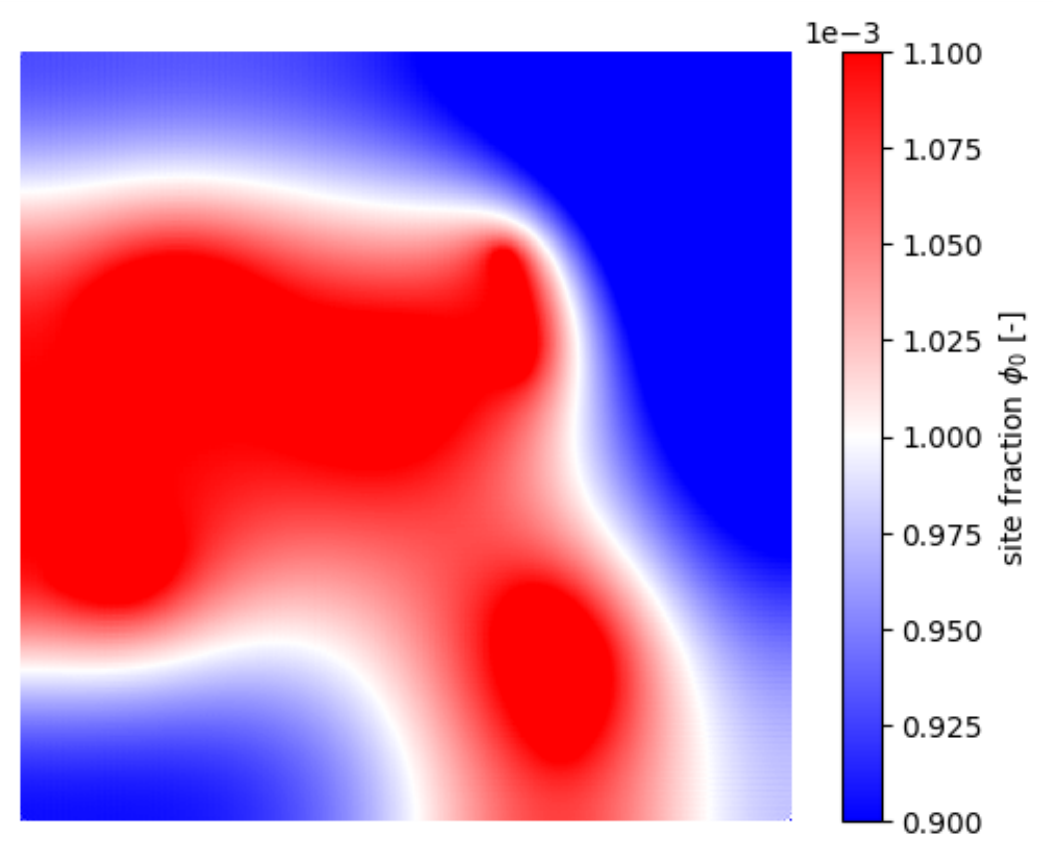}
    \caption{System evolution for the variational problem defined
in equation \eqref{eq:lagrangian_ternary}. The contour visible is the distribution of the site
fraction of vacancies \(\phi_0\) after \(0.1[s]\). The parameters were
chosen as follows: \(\phi_{0}^{eq}=0.001\),
\(D_0 \rightarrow \infty\), \(D_1=2\), \(D_2=1\), $k_0=200[\tfrac{J}{mol}]$, $k_1 = k_2 = 2[\tfrac{J}{mol}]$, $A_\psi = 1/2$ and \(\kappa = 0\).}
    \label{fig12}
    \end{figure}

The interplay between these mechanisms drives the interface
surface area towards minimization in the multidimensional case, as will
be demonstrated in the two-dimensional example. In this study, we
consider a square domain of unit length, initialized with a prescribed
pattern of circular inclusions for the mole fractions as depicted in Figure~\ref{fig7}. The binary
system is composed of two components that are entirely separated from
each other at this stage. The blue region represents 100\% component 2
and phase $\alpha$, while the red region corresponds to 100\% component 1 and
phase $\beta$. The evolution of the classical Cahn-Hilliard model is shown in
Figure~\ref{fig8}. As expected, the mole fraction approaches its equilibrium
value for each present phase, respectively. Furthermore, the model
exhibits a tendency to minimize the contact area between the two phases,
resulting in typical circular shapes. In contrast, we observe that
limiting cases (i) and (ii) exhibit distinct patterns. In these cases,
there are regions where the interface exhibits higher curvature,
indicating that the model does not prioritize contact area minimization.
Nevertheless, the mole fraction reaches equilibrium values for each
phase respectively, as in the Cahn-Hilliard model. For limiting case
(i), Figure~\ref{fig10}, the overall evolution of the structure corresponds to the evolution
of case (ii), Figure~\ref{fig9}, with the notable exception being the smooth
representation of the interface.
Figure~\ref{fig9b} repeats the computation of Figure~\ref{fig9} with the Raviart--Thomas/discontinuous pair discussed in Section~\ref{sec:spaces}. The two element choices yield the same morphology; the mixed pair resolves the interface within a single cell, the equal-order pair within a few elements, which is the sharpest representation a fixed mesh admits. The costs differ considerably. On a single thread, the equal-order computations of Figures~\ref{fig8}, \ref{fig9} and~\ref{fig10} took 17, 109 and 100 minutes and the ternary computation of Figures~\ref{fig11} and~\ref{fig12} 12 hours, whereas the Raviart--Thomas pair needed 4 minutes for Figure~\ref{fig9b} and 14 minutes for the ternary problem, and comparable times for the regularized cases; the complete timings are listed with the published code. The cost is governed by the convergence of the Newton iteration, which requires damping in the convex-hull cases, rather than by the regularization.

Finally, we consider the ternary system with vacancies, defined by the
variational problem in equation \eqref{eq:lagrangian_ternary}. The same initial conditions of the
previous example are applied for component \(\phi_1\). The initial
condition for the vacancies will be uniformly set to its equilibrium
site fraction \(\phi_0(t=0) = \phi_{0}^{eq} = 10^{-3}[-]\). Note that the dissipation term of
vacancy flux is, unlike suggested in \cite{svoboda2006diffusion},
neglected (implying that \(D_0 \rightarrow \infty\)) because it is
assumed that moving a vacancy by itself does not account for any
dissipation. However, since diffusion in substitutional alloys is
primarily enabled by atoms switching sites with vacancies, we
nevertheless expect dissipation which is accounted for by the other
terms, containing \(D_1\) and \(D_2\).
The coefficients are selected as \(D_1 = 2 \, [\text{mm}^2/\text{s}]\) and \(D_2 = 1 \, [\text{mm}^2/\text{s}]\) to induce the Kirkendall effect. This choice of diffusion coefficients results in component \(\phi_1\) diffusing at a faster rate than component \(\phi_2\), thereby creating a disparity that necessitates counterbalancing by a vacancy flux (according to
equation \eqref{eq:flux_constraint}). As can be seen in Figure~\ref{fig12}, the flux of vacancies
\(\underline{j}_0\) will create an excess in vacancies in the area of
the inclusions and a lack of vacancies in the bulk regions, as expected.
The overall evolution of \(\phi_1\) is depicted in Figure~\ref{fig11} and it can be
seen by comparison to Figure~\ref{fig9} that the overall evolution is not
drastically changed, but rather slowed down by the involved vacancy
mechanism.

\section{Discussion}

\subsection{Operating on the boundary of degeneracy}
\label{sec:discussion_degeneracy}

The computations are consistent with the picture of Section~\ref{sec:degeneracy}. For the convex hull without regularization the mole fraction develops a jump, the affinity stays continuous, and everything that happens inside the two-phase region is dictated by the fluxes that the bulk phases impose at its boundary. The scheme does not resolve the interface as a free boundary; it resolves it within a few elements, and the price for this is the indeterminacy of the mole fraction inside the degenerate region, visible as the fluctuations of Figure~\ref{fig3} and as the slightly irregular contour of Figure~\ref{fig9}. This is the behavior known from the enthalpy formulation of the Stefan problem on a fixed mesh \cite{elliott1987error, nochetto1988approximation}, in which the position of the interface is recovered only to the accuracy of a mesh cell and a mushy region of positive measure is consumed from its boundaries. The alternative, a sharp-interface or moving-mesh method \cite{dreyer2017sharp, baines2009moving, baines2011velocity}, removes the degeneracy by construction, at the cost of a mesh that follows the interface, with remeshing, topology changes and field transfer at every step. For the problems we have in mind, diffusional phase transformations coupled to damage and to continuum mechanics in geometries with many inclusions, this cost is prohibitive, and a fixed mesh is the only practical choice.

What the fixed-mesh description retains is the physics. The scheme conserves mass exactly and satisfies the discrete energy--dissipation inequality of Proposition~\ref{prop:consistency} for every convex molar free energy, the convex hull included: the free energy decreases by at least the energy dissipated in each step, and the balance closes exactly once the numerical dissipation of the backward Euler step is accounted for. Neither statement requires the equation to be parabolic, and neither is affected by the degeneracy at the phase boundary. This is the sense in which we operate on the boundary of degeneracy, and it is sufficient for what we are doing: the quantities a thermodynamic model is accountable for, mass, energy and entropy production, are under control, while the fine structure of the interface is not resolved. A more detailed description of the interface would have to come from atomistic models, which lie outside the thermodynamic, continuum description to which we restrict ourselves. We note that the regularized cases (i) and the Cahn--Hilliard model are not degenerate, and that the scheme treats all three cases with the same code; the choice between them is a modeling decision about the interface energy, not a numerical one.

\subsection{On the choice of the finite element spaces}
\label{sec:spaces}

The variational forms \eqref{eq:var_form_1}--\eqref{eq:var_form_3} leave the finite element spaces open, and two choices are natural. Both satisfy Proposition~\ref{prop:consistency}, and both have been implemented with the same code, the only differences being the element definitions and the facet-based discretization of the gradient energy. The equal-order continuous spaces used for the figures of this work are desirable for the coupling to continuum mechanics in established finite element frameworks such as ABAQUS, which provide continuous elements only, and they give a smooth representation of the two-dimensional interfaces; their stability for the Darcy-type subproblem rests on the time-derivative term, as discussed in Section~\ref{subsec:fem}. The inf-sup stable choice for the dual mixed form is the Raviart--Thomas space for the flux (RTCF on quadrilaterals) together with piecewise constant spaces for the mole fraction and the affinity \cite{raviart1977mixed, boffi2013mixed}, in which the flux boundary condition $\underline{j} \cdot \underline{n} = 0$ is imposed on the normal component only, the mass balance holds element-wise, and the mole fraction is allowed to be discontinuous, which is the regularity the degenerate case actually has. For piecewise constant $\phi_h$ the gradient vanishes inside every element, and a gradient energy, if present, is carried by the jumps $\jump{\phi_h}$ across the interior facets $e \in \mathcal{E}_h$ of the mesh. In the spirit of interior penalty methods \cite{arnold1982interior} it is approximated by the facet term
\begin{equation}
    \label{eq:jump_penalty}
    \int_{V} \frac{\kappa}{2} |\nabla \phi|^2 \; dV \;\approx\; \sum_{e \in \mathcal{E}_h} \frac{\kappa}{2 h_e} \int_{e} \jump{\phi_h}^2 \; dS,
\end{equation}
where $h_e$ is the distance between the centres of the two elements sharing the facet $e$, so that the difference quotient $\jump{\phi_h}/h_e$ replaces the normal derivative; accordingly, the term $\kappa \nabla \phi \cdot \nabla \hat{\phi}$ in \eqref{eq:var_form_1} is replaced by $\sum_{e} \frac{\kappa}{h_e} \int_e \jump{\phi}\, \jump{\hat{\phi}} \; dS$. The facet term is a convex functional of $\phi_h$ and part of the discrete free energy, so Proposition~\ref{prop:consistency} holds verbatim with $\mathcal{F}$ replaced by its discrete counterpart; this is the form in which the identity of Proposition~\ref{prop:consistency}(b) was verified numerically for the mixed pair. For $\kappa = 0$, the case of Figure~\ref{fig9b} and Table~\ref{tab:elements}, no facet term enters at all. We regard this choice as the future of the method: it is conservative in the strict, element-wise sense, its flux is a genuine $H(\mathrm{div})$ field, and it is the natural starting point for a moving-mesh extension, in which the mass balance across a moving facet is again a flux jump condition; the coupling of a Raviart--Thomas flux to a continuous displacement field is standard in poroelasticity \cite{phillips2007coupling}. Table~\ref{tab:elements} compares the two choices, and two further variants, on the one-dimensional problem of Figure~\ref{fig3}. The mixed pair gives the sharpest and, at the end, perfectly monotone interface, but the largest fluctuations during the transient: without any gradient control of the piecewise constant mole fraction, the indeterminacy of the degenerate region is fully exposed, whereas the continuity of the equal-order spaces damps it. Figure~\ref{fig9b} shows that the two choices agree on the morphology in two dimensions, at a fraction of the cost for the mixed pair.

\begin{table}[H]
\centering
\caption{Element choices for the one-dimensional case (ii) of Figure~\ref{fig3} (200 cells, $\Delta t = 10^{-3}$, 200 steps). Oscillation: total variation of $\phi_h$ minus its range, which vanishes for a monotone profile; maximum over the transient and value at $t = 0.2$. Width: number of degrees of freedom with $0.26 < \phi_h < 0.74$ at $t = 0.2$. In one dimension the Raviart--Thomas space of order $k$ is the continuous $P_{k+1}$ space.}
\label{tab:elements}
\begin{tabular}{lcccc}
\hline
spaces for $(\phi, \underline{j}, \bar{\mu})$ & max.\ oscillation & final oscillation & width & Newton it./step \\
\hline
$P_1$ / $P_1$ / $P_1$ (equal order) & 0.36 & $4 \cdot 10^{-5}$ & 1 & 1.3 \\
$P_2$ / $P_2$ / $P_2$ (equal order, this work) & 0.50 & $7 \cdot 10^{-3}$ & 5 & 1.8 \\
$P_0^{\mathrm{dc}}$ / $P_1$ / $P_0^{\mathrm{dc}}$ (Raviart--Thomas) & 11.9 & 0 & 0 & 1.9 \\
$P_1^{\mathrm{dc}}$ / $P_2$ / $P_1^{\mathrm{dc}}$ (Raviart--Thomas, next order) & 1.14 & $4 \cdot 10^{-2}$ & 2 & 1.5 \\
\hline
\end{tabular}
\end{table}

\section{Conclusion}

In this paper, we have introduced a novel methodology for mathematically
modeling diffusional phase transformations with and without interface energy. By relating our
approach to limiting cases of the Cahn-Hilliard model, we demonstrate
the versatility and potential of our method in handling complex phase
transformations. The method is derived from the thermodynamic extremal principle: the evolution equations are the stationarity conditions of the Rayleighian, the rate of the free energy plus one half of the dissipation, with respect to the rates, fluxes and multipliers at a frozen state, and the time-discrete scheme is the corresponding incremental minimization. The incorporation and explicit solving for
rate-dependent variables like the flux allows us to avoid differentiating the affinity
and to overcome limitations imposed by traditional numerical techniques. We have made explicit that the convex hull without regularization, the CALPHAD limit, is a degenerate parabolic problem of Stefan type, that the degeneracy can be removed only by a moving mesh, which is incompatible with our aim of coupling diffusional phase transformations and damage to continuum mechanics, and that operating on the boundary of degeneracy is nevertheless sound: the scheme conserves mass exactly and satisfies a discrete energy--dissipation inequality, the second law, for every convex free energy and for every choice of conforming finite element spaces.

In our opinion this degeneracy has to be taken care of numerically, but it is not a physical problem, because the continuity equation $\frac{\partial \phi}{\partial t} + \nabla \cdot \underline{j} = 0$ does not have to be satisfied in the strong, pointwise sense in the first place. Thermodynamics presupposes a large number of particles: the mole fraction, the affinity and the flux are averages over volume elements that contain many atoms, and in describing a microstructure in all its details one reaches a barrier below which these fields, and with them their derivatives, lose their meaning. What conservation of mass and thermodynamics can assert are balances for such volume elements: the matter that leaves one of them enters its neighbours, and the total free energy does not increase. The strong form follows from these balances only for a continuum in the abstract sense, whose fields can be differentiated. This is too strong a requirement as soon as particles are involved, and the continuum model itself violates it at a sharp interface, where the mole fraction jumps and neither its rate nor the divergence of the flux exists as a function. The weak form, \eqref{eq:stat_3} and \eqref{eq:var_form_3}, is sufficient. It demands the balance only for averages over finite volumes, and the jump condition that a sharp-interface description imposes at the interface is contained in it. With the piecewise constant test functions of the mixed pair of Section~\ref{sec:spaces} it is literally the balance of every mesh cell with the fluxes through its faces; with the continuous spaces it holds for the weighted averages defined by the test functions and, by Proposition~\ref{prop:consistency}, for the body as a whole. This is also the position of the mathematical theory, in which the solutions of problems of Stefan type are defined in the weak sense \cite{visintin1996models, elliott1987error}. The weak form is therefore not a relaxation of the physical statement but the physical statement itself.

Our approach to diffusional phase transformations and the phase-field
method provides a novel finite element
approach, offering a promising framework for investigating these
phenomena. We also demonstrated the thermodynamic consistency of the
reactive diffusion models and the phase-field method in general and
helped to bridge the gaps between these fields. In doing so we also
removed a widespread misconception regarding the chemical potential and
its role as driving force for diffusion. Our work contributes to the
development of new mathematical tools for modeling complex phase
transformations in materials science, ultimately enriching our
understanding of these fundamental processes. Future studies can build
upon this foundation, exploring further applications and refinements to
enhance the predictive capabilities of our methodology.

\subsection{CRediT authorship contribution statement}

\textbf{W. Flachberger:} Conceptualization, Methodology, Software, Validation, Formal analysis, Investigation, Data curation, Visualization, Writing -- original draft, Writing -- review \& editing (lead). \textbf{T. Antretter, J. Svoboda, S. Leitner, M. Petersmann:} Supervision, Project administration, Funding acquisition, Writing -- review \& editing (supporting). \textbf{S. Gaddikere-Nagaraja:} Writing -- review \& editing (supporting).

\subsection{Acknowledgement}

    The authors gratefully acknowledge the financial support under the scope
of the COMET program within the K2 Center ``Integrated Computational
Material, Process and Product Engineering (IC-MPPE)'' (Project 886385).
This program is supported by the Austrian Federal Ministries for Climate
Action, Environment, Energy, Mobility, Innovation and Technology (BMK)
and for Digital and Economic Affairs (BMDW), represented by the Austrian
research funding association (FFG), and the federal states of Styria,
Upper Austria and Tyrol.

Additionally the authors thank Jeremy Bleyer,
currently affiliated with Laboratoire Navier of the Ecole des Ponts
ParisTech, who gave valuable advice on mixed finite element
formulations.

\appendix
\section{Appendix: Demonstrating the Satisfaction of the Gibbs-Duhem Condition}
\label{app:Gibbs-Duhem}
    The Gibbs-Duhem condition reads:

\begin{equation}
    \label{eq:app_gibbs-duhem}
    \sum_{j=1}^{n} m_j  d\mu_j = 0
\end{equation}

For a mass conserved system, equation \eqref{eq:app_gibbs-duhem} can be divided by the total
number of moles to yield an expression including the mole fraction
(according to equations \eqref{eq:total_number_moles} and \eqref{eq:mole_fraction}):

\begin{equation}
    \sum_{j=1}^{n} \phi_j  d\mu_j = 0.
\end{equation}

Again, the last element of the sum can be explicitly addressed:

\begin{equation}
    \sum_{j=1}^{n-1} \phi_j  d\mu_j  + \phi_n  d\mu_n = 0.
\end{equation}

The chemical potential from equation \eqref{eq:mu_i-x_i} can be inserted for \(\mu\):

\begin{equation}
    \sum_{j=1}^{n-1} \phi_j  d\Big( f_m + \sum_{k=1}^{n-1} (\delta_{jk} - \phi_k) \frac{\partial f_m}{\partial \phi_k} \Big)  +  \phi_n  d\Big(  f_m - \sum_{k=1}^{n-1} \phi_k \frac{\partial f_m}{\partial \phi_k} \Big) = 0
\end{equation}

The mole fraction of the last component was defined in equation \eqref{eq:x_n} and
can be substituted for \(\phi_n\):

\begin{equation}
    \sum_{j=1}^{n-1} \phi_j \Big( df_m + \sum_{k=1}^{n-1} (\delta_{jk} - \phi_k)  d\Big(\frac{\partial f_m}{\partial \phi_k}\Big) - \sum_{k=1}^{n-1} d\phi_k  \Big(\frac{\partial f_m}{\partial \phi_k}\Big) \Big) + \big( 1 - \sum_{i=1}^{n-1} \phi_i \big)  \Big(  df_m - \sum_{k=1}^{n-1} d\phi_k \frac{\partial f_m}{\partial \phi_k} - \sum_{k=1}^{n-1} \phi_k d\Big(\frac{\partial f_m}{\partial \phi_k} \Big) \Big) = 0
\end{equation}

\begin{equation}
    \sum_{j=1}^{n-1} \phi_j \Big( \sum_{k=1}^{n-1} (\delta_{jk} - \phi_k)  d\Big(\frac{\partial f_m}{\partial \phi_k}\Big)  \Big) + \big( 1 - \sum_{i=1}^{n-1} \phi_i \big)  \Big(  - \sum_{k=1}^{n-1} \phi_k d\Big(\frac{\partial f_m}{\partial \phi_k} \Big) \Big)  = 0
\end{equation}

\begin{equation}
    \sum_{j=1}^{n-1} \sum_{k=1}^{n-1} \phi_j (\delta_{jk} - \phi_k)  d\Big(\frac{\partial f_m}{\partial \phi_k}\Big) - \sum_{k=1}^{n-1} \phi_k d\Big(\frac{\partial f_m}{\partial \phi_k} \Big) + \sum_{i=1}^{n-1} \sum_{k=1}^{n-1} \phi_i \phi_k d\Big(\frac{\partial f_m}{\partial \phi_k} \Big)  = 0
\end{equation}

\begin{equation}
    \sum_{j=1}^{n-1} \sum_{k=1}^{n-1} \phi_j \delta_{jk} d\Big(\frac{\partial f_m}{\partial \phi_k}\Big) - \sum_{j=1}^{n-1} \sum_{k=1}^{n-1} \phi_j \phi_k d\Big(\frac{\partial f_m}{\partial \phi_k}\Big) - \sum_{k=1}^{n-1} \phi_k d\Big(\frac{\partial f_m}{\partial \phi_k} \Big) + \sum_{i=1}^{n-1} \sum_{k=1}^{n-1} \phi_i \phi_k d\Big(\frac{\partial f_m}{\partial \phi_k} \Big)  = 0
\end{equation}

\begin{equation}
    \sum_{j=1}^{n-1} \sum_{k=1}^{n-1} \phi_j \delta_{jk} d\Big(\frac{\partial f_m}{\partial \phi_k}\Big) - \sum_{k=1}^{n-1} \phi_k d\Big(\frac{\partial f_m}{\partial \phi_k} \Big)= 0
\end{equation}

\begin{equation}
    \sum_{j=1}^{n-1}  \phi_j d\Big(\frac{\partial f_m}{\partial \phi_j}\Big) - \sum_{k=1}^{n-1} \phi_k d\Big(\frac{\partial f_m}{\partial \phi_k} \Big)= 0
\end{equation}

\qed

\section{Appendix: Analytical expressions for the molar free energies \texorpdfstring{$f_0$}{f0} and \texorpdfstring{$f_{0}^{**}$}{f0**}}
\label{app:functions}

\begin{equation}
    f_{0}^{**}(\phi) = \tfrac{2}{100} \phi + (-\phi + ((\tfrac{1}{4} - \phi)^2)^{1/2} + \tfrac{1}{4})^2/4 + (\phi + ((\phi - \tfrac{3}{4})^2)^{1/2} - \tfrac{3}{4})^2/4 + \tfrac{5}{1000}
\end{equation}

\begin{multline}
    f_{0}(\phi) = \tfrac{2}{100}\phi + (-\phi + ((\tfrac{1}{4} - \phi)^2)^{1/2} + \tfrac{1}{4})^2/4 \\
    + (\phi + ((\phi - \tfrac{3}{4})^2)^{1/2} - \tfrac{3}{4})^2/4 \\
    + (-(\tfrac{3}{4} - \phi)/(2((\phi - \tfrac{3}{4})^2)^{1/2}) \\
    - (2\phi - \tfrac{3}{2})/(2((\phi - \tfrac{3}{4})^2)^{1/2}) \\
    - (-16(\phi - \tfrac{1}{4})^3 \\
    + 12(\phi - \tfrac{1}{4})^2)((\tfrac{1}{4} - \phi)/(2((\phi - \tfrac{1}{4})^2)^{1/2}) \\
    - (\tfrac{3}{4} - \phi)/(2((\phi - \tfrac{3}{4})^2)^{1/2}) \\
    - (2\phi - \tfrac{3}{2})/(2((\phi - \tfrac{3}{4})^2)^{1/2}) \\
    + (2\phi - \tfrac{1}{2})/(2((\phi - \tfrac{1}{4})^2)^{1/2})) \\
    + \tfrac{1}{2})((\tfrac{3}{4} - \phi)/(2((\phi - \tfrac{3}{4})^2)^{1/2}) \\
    + (2\phi - \tfrac{3}{2})/(2((\phi - \tfrac{3}{4})^2)^{1/2}) \\
    + (-16(\phi - \tfrac{1}{4})^3 \\
    + 12(\phi - \tfrac{1}{4})^2)((\tfrac{1}{4} - \phi)/(2((\phi - \tfrac{1}{4})^2)^{1/2}) \\
    - (\tfrac{3}{4} - \phi)/(2((\phi - \tfrac{3}{4})^2)^{1/2}) \\
    - (2\phi - \tfrac{3}{2})/(2((\phi - \tfrac{3}{4})^2)^{1/2}) \\
    + (2\phi - \tfrac{1}{2})/(2((\phi - \tfrac{1}{4})^2)^{1/2})) + \tfrac{1}{2})/30 + \tfrac{5}{1000}
\end{multline}

\bibliographystyle{unsrtnat}
\bibliography{references}

\end{document}